\documentclass[reqno, 11pt]{amsart}  
\usepackage{amsmath}
\usepackage{amssymb}
\usepackage{hyperref}
\hypersetup{colorlinks=true}

\usepackage{graphicx}

\usepackage[left=1.125in,right=1.1in,top=1.1in,bottom=1.125in]{geometry}

\newtheorem{prop}{Proposition}[section]
\newtheorem{lemma}[prop]{Lemma}
\newtheorem{thm}[prop]{Theorem}
\newtheorem{cor}[prop]{Corollary}
\numberwithin{equation}{section}

\newcommand{\shrink}[1]{ {\scriptstyle {\textstyle {#1} } } }
\newcommand{\smfrac}[2]{ \shrink{ \frac{#1}{#2} } }
\newcommand{\lin}{\langle} 
\newcommand{\rin}{\rangle} 

\newcommand{\nt}{\negthinspace}

\newcommand{\fom}{\mathtt{fom}} 
\newcommand{\parfom}{\parallel \nt \nt  \mathtt{FOM}}
\newcommand{\subgrad}{\mathtt{subgrad}} 
\newcommand{\accel}{\mathtt{accel}} 
\newcommand{\fista}{\mathtt{fista}} 
\newcommand{\univ}{\mathtt{univ}}  
\newcommand{\smooth}{\mathtt{smooth}}

\newcommand{\tpause}{t_{\mathrm{pause}}} 
\newcommand{\ttransit}{t_{\mathrm{transit}}}

\newcommand{\sparfom}{\mathtt{Sync} \nt \nt  \parfom} 
\newcommand{\aparfom}{\mathtt{Async} \nt \nt \parfom}

\newcommand{\grad}{\nabla} 
 
\newcommand{\dist}{\mathrm{dist}} 

\newcommand{\Sym}{ \mathbb{S}^n } 

\usepackage{array,multirow}

\begin{document}

\newpage 
$ \textrm{~} $ \quad \vspace{-3mm}

\title[A Simple Nearly-Optimal Restart Scheme]{A Simple Nearly-Optimal Restart Scheme \\ For Speeding-Up First Order Methods}

\begin{abstract}
We present a simple scheme for restarting first-order methods for convex optimization problems. Restarts are made based only on achieving specified decreases in objective values, the specified amounts being the same for all optimization problems. Unlike existing restart schemes, the scheme makes no attempt to learn parameter values characterizing the structure of an optimization problem, nor does it require any special information that would not be available in practice (unless the first-order method chosen to be employed in the scheme itself requires special information).  As immediate corollaries to the main theorems, we show that when some well-known first-order methods are employed in the scheme, the resulting complexity bounds are  nearly optimal for particular -- yet quite general -- classes of problems.
\end{abstract} \vspace{-3mm}

\author[J. Renegar and B. Grimmer]{James Renegar and Benjamin Grimmer}
\address{School of Operations Research and Information Engineering,
 Cornell University, Ithaca, NY, U.S.} 
\thanks{Research supported in part by NSF grants CCF-1552518 and DMS-1812904, and by the NSF Graduate Research Fellowship under grant DGE-1650441.
A portion of the initial development was done while the authors were visiting the Simons Institute for the Theory of Computing, and was partially supported by the DIMACS/Simons Collaboration on Bridging Continuous and Discrete Optimization through NSF grant CCF-1740425. }

\thanks{We are indebted to reviewers, whose comments and suggestions led to a wide range of substantial improvements in the presentation.}

\maketitle

\vspace{-7mm}

\section{{\bf  Introduction}}  \label{sect.a}

A restart scheme is a procedure that stops an algorithm when a given criterion is satisfied, and then restarts the algorithm with new input.  For certain classes of convex optimization problems, restarting can improve the convergence rate of first-order methods, an understanding that goes back decades to work of Nemirovski and Nesterov \cite{NemNes85}. Their focus, however, was on the abstract setting in which somehow known are various scalars from which can be deduced nearly-ideal times to make restarts for the particular optimization problem being solved. 

In recent years, the notion of ``adaptivity'' has come to play a foremost role in research on first-order methods. Here a scheme learns -- when given an optimization problem to solve -- good times at which to restart the first-order method by systematically trying various values for the methods's parameters and observing consequences over a number of iterations.  Such a scheme, in the literature to date, is particular to a class of optimization problems, and typically is particular to the first-order method to be restarted. No simple and easily-implementable set of principles has been developed which applies to arrays of first-order methods and arrays of problem classes. We eliminate this shortcoming of the literature by introducing such a set of principles.
 
We present a simple restart scheme utilizing multiple instances of a general first-order method. The instances are run in parallel (or sequentially if preferred) and occasionally communicate their improvements in objective value to one another, possibly triggering restarts. In particular, a restart is triggered only by improvement in objective value, the trigger threshold being independent of problem classes and first-order methods. Nevertheless, for a variety of prominent first-order methods and important problem classes, we show nearly-optimal complexity is achieved.

We begin by outlining the key ideas, and providing an overview of ramifications for theory.  

\subsection{The setting}  \label{sect.aa}

We consider optimization problems
\begin{equation}  \label{eqn.aa} 
\begin{array}{rl}
\min & f(x) \\
\textrm{s.t.} & x \in Q \; , \end{array} \end{equation}
where $ f $ is a convex function, and $ Q $ is a closed convex set contained in the (effective) domain of $ f $ (i.e., where $ f $ is finite). We assume the set of optimal solutions, $ X^* $, is nonempty. Let $ f^* $ denote the optimal value.

Let $ \fom $ denote a first-order method capable of solving some class of convex optimization problems of the form (\ref{eqn.aa}), in the sense that for any problem in the class, when given an initial feasible point $ x_0 \in Q $ and accuracy $ \epsilon  > 0 $, the method is guaranteed to generate a feasible iterate $ x_k $ satisfying $ f(x_k) \leq f^* + \epsilon $, an ``$ \epsilon $-optimal solution.'' 

For example, consider the class of problems for which $ f $ is Lipschitz on an open neighborhood of $ Q $, that is, there exists $ M > 0 $ such that $ | f(x) - f(y) | \leq M \| x - y \| $ for all $ x $, $ y $ in the neighborhood (where, throughout, $ \| \; \; \| $ is the Euclidean norm). An appropriate algorithm is the projected subgradient method, $ \fom = \subgrad $, defined iteratively by
\begin{equation}  \label{eqn.ab} 
   x_{k+1} = P_Q \big( x_k - \smfrac{\epsilon}{\| g_k \|^2} \, g_k \big) \; , 
   \end{equation} 
where $ g_k $ is any subgradient\footnote{Recall that a vector $ g $ is a subgradient of $ f $ at $ x $ if for all $ y $ we have $ f(y) \geq f(x) + \lin g, y-x \rin $.  If $ f $ is differentiable at $ x $, then $ g $ is precisely the gradient.} of $ f $ at $ x_k $, and where $ P_Q $ denotes orthogonal projection onto $ Q $, i.e., $ P_Q(x) $ is the point in $ Q $ which is closest to $ x $. Here the algorithm depends on $ \epsilon  $ but not on the Lipschitz constant, $ M $.

As another example, consider the class of problems for which $ f $ is $ L $-smooth on an open neighborhood of $ Q $, that is, $ f $ is differentiable and satisfies $ \| \grad f(x) -\grad f(y) \| \leq L \| x - y \| $ for all $ x $, $ y $ in the neighborhood. These are the problems for which accelerated methods were initially developed, beginning with the pioneering work of Nesterov \cite{nesterov1983method}. While our framework is immediately applicable to many accelerated methods, for concreteness we refer to a particular instance of the popular algorithm FISTA \cite{beck2009fast}, a direct descendent of Nesterov's original accelerated method. We signify this instance as $ \accel $ and record it below\footnote{$ \accel $: \\ $ \textrm{~} $ \qquad (0)  Input: $ x_0 \in Q $. Initialize: $ y_0 = x_0 $, $ \theta_0 = 1 $, $ k = 0 $. \\
$ \textrm{~} $ \qquad (1) $ x_{k+1} = P_Q( y_k - \smfrac{1}{L}  \grad f(y_k)) $ \\
$ \textrm{~} $ \qquad (2) $ \theta_{k+1} = \smfrac{1}{2} (1 + \sqrt{1 + 4 \theta_k^2}) $ \\
$ \textrm{~} $ \qquad (3) $ y_{k+1} = x_{k+1} + \smfrac{\theta_k - 1}{\theta_{k+1}} (x_{k+1} - x_k) $ }. We note that $ \accel $  depends on $ L $ but not on $ \epsilon $. (Later we discuss ``backtracking,'' where only an estimate $ \bar{L} $ of $ L $ is required, with an algorithm automatically updating $ \bar{L} $ during iterations, increasing the overall complexity by only an additive term depending on the quantity $ | \log( \bar{L}/L ) | $.)  

The basic information required by $ \subgrad $ and $ \accel $ are (sub)gradients. Both methods make one call to the ``oracle'' per iteration, to obtain that information for the current iterate. 

Our theorems are general, applying to any first-order method for which there exists an upper bound on the total number of oracle calls sufficient to compute an $ \epsilon $-optimal solution, where the upper bound is a function of $ \epsilon $ and the distance from the initial iterate $ x_0 $ to the set of optimal solutions.  Moreover, the theorems allow any oracle. 

Thus, for example, besides immediately yielding corollaries specific to $ \subgrad $ and $ \accel $, our theorems readily yield corollaries for accelerated proximal-gradient methods, which are relevant when $ Q = \mathbb{R}^n $, $ f = f_1 + f_2 $ with $ f_1 $ being convex and smooth, and $ f_2 $ being convex and ``nice'' (e.g., a regularizer). Here the oracle provides an optimal solution to a particular subproblem involving a gradient of $ f_1 $ and the ``proximal operator'' for $ f_2 $ (c.f., \cite{nesterov2013gradient}). For example, our results apply to $ \mathtt{fista} $, the general version of $ \accel $, as will be clear to readers familiar with proximal methods.\footnote{We avoid digressing to discuss the principles underlying proximal methods.  For readers unfamiliar with the topic, the paper introducing FISTA (\cite{beck2009fast}) is a good place to start.} 

\subsection{Brief motivation for the scheme}  \label{sect.ab} 

To motivate our restart scheme, consider the subgradient method, $ \fom = \subgrad $. If $ \epsilon $ is small and $ f(x_0) - f^* $ is large, the subgradient method tends to be extremely slow in reaching an $ \epsilon $-optimal solution, in part due to $ x_0 $ being far from $ X^* $ and the step sizes $ \| x_{k+1} - x_k \|  $ being proportional to $ \epsilon $.  In this situation it would be preferable to begin by relying on a constant $ 2^N \epsilon $  in place of $ \epsilon $, where $ N $ is a large integer, continuing to make iterations until a $ (2^N \epsilon) $-optimal solution has been reached, then rely on the value $ 2^{N-1} \epsilon   $ instead, iterating  until a $ (2^{N-1} \epsilon)  $-optimal solution has been reached, and so on, until finally a $ (2^0 \epsilon) $-optimal solution has been reached, at which time the procedure would terminate.

The catch, of course, is that generally there is no efficient way to determine whether a $ (2^n \epsilon) $-optimal solution has been reached, no efficient termination criterion. To get around this conundrum, our scheme relies on the value $ f(x_0) - f(x_k) $, that is, relies on the  decrease in objective {\em  that has already been achieved}, and does not attempt to measure how much objective decrease remains to be made. 

\subsection{The key ideas}  \label{sect.ac}

For $ n = -1, 0, 1, \ldots, N $, let $ \fom_n $ be a version of $ \fom $ which is known to compute an iterate which is a $ (2^n \epsilon) $-optimal solution. In particular, when $ \fom = \subgrad $, let $ \fom_n $ be the subgradient method with $ 2^n \epsilon $ in place of $ \epsilon $, and when $ \fom = \accel $, simply let $ \fom_n = \accel $ (the accelerated method is independent of $ \epsilon $).

For definiteness, assume the desired accuracy $ \epsilon $ satisfies $ 0 < \epsilon \ll 1 $, and choose $ N $ so that $ 2^N \approx 1/\epsilon $. Thus, $ N \approx \log_2(1/\epsilon) $.

The algorithms $ \fom_n $ are run in parallel (or sequentially if preferred), with $ \fom_n $ being initiated at a point $ x_{n,0} $. Each algorithm is assigned a task. The task of $ \fom_n $ is to make progress $ 2^n \epsilon $, that is, to obtain a feasible point $ x $  satisfying $ f(x) \leq f(x_{n,0}) -  2^n \epsilon $. If the task is accomplished, $ \fom_n $ notifies $ \fom_{n-1} $ of the point $ x $ (assuming $ n \neq -1 $), to allow for the possibility that $ x $ also serves to accomplish the task of $ \fom_{n-1} $. 

If $ n \neq N $, then in addition to $ \fom_{n-1} $ being notified of $ x $, a restart for $ \fom_n $ occurs at $ x $, that is, $ x $ becomes the (new) initial point $ x_{n,0} $ for $ \fom_n $. Again the assigned task of $ \fom_n $ is to obtain a point $ x $ satisfying $ f(x) \leq f(x_{n,0}) -  2^n \epsilon $. (The algorithm $ \fom_N $ never restarts, for reasons particular to obtaining complexity bounds of a particular form.)
\vspace{2mm}

\begin{quote} 
Thus, in a nutshell, the main ideas for the scheme are simply that (1) $ \fom_n $ continues iterating either until reaching an iterate $ x = x_{n,k} $ satisfying the inequality $ f(x) \leq f(x_{n,0}) -  2^n \epsilon $ (where $ x_{n,0} $ is the most recent restart point for $ \fom_n $), or until receiving a point $ x $ from $ \fom_{n+1} $ which satisfies the inequality, and (2) when $ \fom_n $ obtains such a point $ x $ in either way, then $ \fom_n $ restarts\footnote{In the case of $ \subgrad_n $, restarting at $ x $ simply means computing a subgradient step at $ x $, i.e., $ x_+ = x - \frac{2^n \epsilon }{\| g_k \|^2} g_k $. For $ \accel_n $, restarting is slightly more involved, in part due to two sequences of points being computed, $ \{ x_{n,k} \} $ and $ \{  y_{n,k} \}  $. The primary sequence is $ \{ x_{n,k} \} $. The method $ \accel_n $ restarts at a point $ x $ from a primary sequence (either from $ \fom_n $'s own primary sequence, or a point sent to $ \fom_n $ from the primary sequence for $ \fom_{n+1} $). In restarting, $ \accel_n $ begins by ``reinitializing,'' setting $ x_{n,0} = x $, $ y_{n,0} = x $, $ \theta_0 = 1 $, $ k = 0 $, then begins iterating.} at $ x $ (i.e., $ x_{n,0} \leftarrow x $), and notifies $ \fom_{n-1} $ of $ x $.
\end{quote}
\vspace{2mm}

There are a variety of ways the simple ideas can be implemented, including sequentially, or in parallel, in which case either synchronously or asynchronously. In order to provide rigorous proofs of the scheme's performance, a number of secondary details have to be specified for the scheme, although these details could be chosen in multiple ways and still similar theorems would result. Our choices for the details of the scheme when implemented sequentially, or in parallel with iterations sychronized, are given in \S\ref{sect.d}. Details are given for a parallel asynchronous implementation in \S\ref{sect.g}. The key ingredients to all settings, however, are precisely the simple ideas posited above.

The simplicity of the underlying ideas gives hope for practice. In this regard we provide, in \S\ref{sect.i}, numerical results displaying notable performance improvements when the scheme is employed with either $\mathtt{subgrad}$ or $\mathtt{accel}$. 

\subsection{Consequences for theory} \label{sect.ad}

To theoretically demonstrate effectiveness of the scheme, we consider objective functions $ f $ possessing ``H\"{o}lderian growth,'' in that there exist constants $ \mu > 0 $ and $ d \geq 1 $ for which
\[ 
   x \in Q \textrm{ and } f(x) \leq f( \mathbf{x_0})  \quad \Rightarrow \quad f(x) - f^* \geq \mu \,  \dist(x,X^*)^d \; , 
\] 
where $ \mathbf{x_0} $ is a common point at which each $ \fom_n $ is first started, and where $ \dist(x,X^*) := \min \{ \| x - x^* \| \mid x^* \in X^* \} $. (Elsewhere the property is referred to as a H\"{o}lderian error bound \cite{yang2017adaptive}, as  sharpness \cite{roulet2017sharpness}, as the {\L}ojasiewicz property \cite{karimi2016linear}, and the (generalized) {\L}ojasiewicz inequality  \cite{bolte2007lojasiewicz}.)  We call $ d $ the ``degree of growth.'' The parameters $ \mu $ and $ d $ are not assumed to be known to our scheme, nor is any attempt made to learn approximations to their values. Indeed, our theorems do not assume $ f $ possesses H\"{o}lderian growth.

{\L}ojasiewicz \cite{lojasiewicz1963propriete, lojasiewicz1993geometrie} showed H\"{o}lderian growth holds for generic analytic and subanalytic functions, a result which was generalized to nonsmooth subanalytic convex functions by Bolte, Daniilidis and Lewis \cite{bolte2007lojasiewicz}.  It is easily seen that strongly convex functions have H\"{o}lderian growth with $ d = 2 $, but so do many other (smooth and nonsmooth) convex functions, such as the objective in the Lagrangian form of Lasso regression. Another example to keep in mind is when $ f $ is a piecewise-linear convex function and $ Q $ is polyhedral, then H\"{o}lderian growth holds with $ d = 1 $.

As a simple consequence (Corollary~\ref{cor.eb}) to the theorem for our synchronous scheme (Theorem~\ref{thm.ea}), we show that if $ f $ is $ M $-Lipschitz on an open neighborhood of $ Q $, and $ \subgrad $  is employed in the scheme, the time sufficient to compute an $ \epsilon $-optimal solution (feasible $ x $ satisfying $ f(x) - f^* \leq \epsilon < 1 $) is at most on the order of $ (M/ \mu)^2 \log(1/ \epsilon) $ if $ d = 1 $, and at most order of $ M^2/(\mu^{2/d} \epsilon^{2(1 - 1/d)}) $ if $ d > 1 $. These time bounds apply when the scheme is run in parallel, relying on $ 2 + \lceil \log_2( 1/ \epsilon) \rceil  $ copies of the subgradient method. 

 As there are $ O(\log(1/\epsilon)) $  copies of the subgradient method being run simultaneously, the total amount of work (total number of oracle calls) is obtained by multiplying the above time bounds by $ O(\log(1/\epsilon)) $. 
 
For general convex optimization problems with Lipschitz objective functions, this is the first algorithm which requires only that subgradients are computable (in particular, does not require any information about $ f $ as input, such as the optimal value $ f^* $), {\em  and} which has  the property that if  $ f $ happens to possess linear growth (i.e., happens to be H\"{o}lderian with $ d = 1 $), an $ \epsilon $-optimal solution is computed within a total number of subgradient evaluations depending only logarithmically on $ 1/\epsilon  $. (We review the related literature momentarily.)
 
Another simple consequence of the theorem pertains to the setting where $ f $ is $ L $-smooth on an open neighborhood of $ Q $. (For smooth functions, the growth degree necessarily satisfies $ d \geq 2 $.)  Here we show that if $ \accel $  is employed in the synchronous scheme, the time required to compute an $ \epsilon $-optimal solution is at most on the order of $ \sqrt{L/\mu} \log(1/ \epsilon) $ if $ d = 2 $, and at most on the order of $ L^{\frac{1}{2}}/ (\mu^{\frac{1}{d}} \epsilon^{\frac{1}{2} - \frac{1}{d}}) $    if $ d > 2 $. For the total amount of work, simply multiply by $ O(\log(1/\epsilon)) $. (We also observe the same bounds apply to $ \fista $ in the more general setting of proximal methods.)

In the case $ d = 2 $, our bound on the total number of gradient evaluations (oracle calls) is off by a factor of  only $ O(\log(1/ \epsilon )) $ from the well-known lower bound for strongly-convex smooth functions (\cite{nemirovsky1983problem}), a small subset of the functions having H\"{o}lderian growth with $ d = 2 $.  For $ d > 2 $, the bounds also are off by a factor of only $ O(\log(1/ \epsilon )) $, according to the lower bounds stated by Nemirovski and Nesterov \cite[page 26]{NemNes85}  (for which we know of no proofs recorded in the literature).

A third straightforward corollary to the theorem for our synchronous parallel scheme applies to an algorithm which, following Nesterov\cite{nesterov2005smooth}, makes use of an accelerated method in solving a ``smoothing'' of $ f $, assuming $ f $ to be Lipschitz. We denote the algorithm by $ \smooth $. We avoid digressing to introduce the necessary notation until \S\ref{sect.b}, but here we mention a consequence when employing $ \smooth $ in our synchronous scheme, the goal being to minimize a piecewise-linear convex function
\[ 
      f(x) = \max \{ a_i^T x + b_i \mid i = 1, \ldots, m \}.
\]     
Note that $ M = \max_i \| a_i \| $ is a Lipschitz constant for $ f $. This function has H\"{o}lderian growth with $ d = 1 $ (i.e., linear growth), for some $ \mu > 0 $. In the context of $ M $-Lipschitz convex functions with linear growth, the ratio $ M/\mu $ is sometimes called the ``condition number,'' and is easily shown to have value at least 1.

From our previous discussion, when employing $ \subgrad $ in our synchronous scheme to minimize $ f $, the time bound is of order $ (M/\mu)^2 \log(1/\epsilon) $, quadratic in the condition number. By contrast, employing $ \smooth $ results in a bound of order $ (M/\mu)  \log(1/\epsilon) \sqrt{ \log(m)}$, which generally is far superior, due to being only linearly dependent on the condition number. (The difference is observed empirically in \S\ref{sect.i}.)  

While in theory, smoothings exist for general Lipschitz convex functions, the set of functions known to have {\em  tractable} smoothings is limited. Thus, in practice, $ \subgrad $ is relevant to a wider class of problems than $ \smooth $. Nonetheless, when a smoothing is available, employing $ \smooth $ in our synchronous scheme can be much preferred to employing $ \subgrad $. 
 
Each iteration of $ \subgrad $  requires the same amount of work, one call to the oracle. Likewise for $ \accel $ and $ \smooth $. By contrast, the work required by adaptive methods -- such as methods employing ``backtracking'' -- typically varies among iterations, an especially important example being Nesterov's universal fast gradient method \cite{nesterov2015universal}. Such algorithms can be ill-suited  for use in a synchronous parallel scheme in which each $ \fom_n $ makes one iteration per time period.  Nonetheless, these algorithms fit well in our asynchronous scheme.  As a straightforward consequence (Corollary~\ref{cor.hb}) to the last of our theorems (Theorem~\ref{thm.ha}), we show combining the asynchronous scheme with the universal fast gradient method results in a nearly-optimal number of gradient evaluations when the scheme is applied to the general class of problems in which $ f $ both has H\"{o}lderian growth and has ``H\"{o}lder continuous gradient with exponent $ 0 < \nu \leq 1 $.'' This is the first algorithm for this especially-general class of problems that both is nearly optimal in the number of oracle calls and does not require information that typically would be unavailable in practice.

\subsection{Positioning within the literature on smooth convex optimization.} \label{sect.ae} As remarked above, an understanding that restarting can speed-up first-order methods goes back to work of Nemirovski and Nesterov \cite{NemNes85}, although their focus was on the abstract setting in which somehow known are various scalars from which can be deduced nearly-ideal times to make restarts for the particular optimization problem being solved. 

In the last decade, adaptivity has been a primary focus in research on optimization algorithms, but far more in the context of adaptive accelerated methods than in the context  of adaptive schemes (meta-algorithms). 

Impetus for the development of adaptive accelerated methods was provided by Nesterov in  \cite{nesterov2013gradient}, where he designed and analyzed an  accelerated method in the context of $ f $ being $ L $-smooth, the new method possessing the same convergence rate as his original (optimal) accelerated method (in particular, $ O(\sqrt{L/\epsilon}) $), but without needing $ L $ as input. (Instead, input meant to approximate $ L $ is needed, and during iterations, the method modifies the input value until reaching a value which approximates $ L $ to within appropriate accuracy.)  Secondary consideration, in \S5.3 of that paper, was given specifically to strongly convex functions, with a proposed grid-search procedure for approximating to within appropriate accuracy the so-called strong convexity parameter (as well as $ L $), leading overall to an adaptive accelerated method possessing up to a logarithmic factor the optimal time bound for the class of smooth and strongly-convex functions, demonstrating how designing an adaptive method aimed at a narrower class of functions can result in a dramatically improved convergence rate ($ O(\log(1/\epsilon)) $ vs. $ O(1/ \sqrt{\epsilon}) $), even without having to know apriori the Lipschitz constant $ L $ or the strong convexity parameter). 

A range of significant papers on adaptive accelerated methods followed, in the setting of $ f $ being smooth and either strongly convex or uniformly convex (c.f., \cite{juditsky2014primal}, \cite{lin2014adaptive}), but also relaxing the strong convexity requirement to assuming  H\"{o}lderian growth with degree $ d = 2 $ (c.f., \cite{necoara2016linear}). 

Far fewer have been the number of papers regarding adaptive schemes, but still, important contributions have been made, in earlier papers which focused largely on  heuristics for the 
  setting of $ f $ being $ L $-smooth and strongly convex (c.f., 
 \cite{o2015adaptive,giselsson2014monotonicity,fercoq2016restarting}), but more recently, papers analyzing schemes with proven guarantees (c.f., \cite{fercoq2016restarting}), including for the setting when the smooth function $ f $ is not necessarily strongly convex but instead satisfies the weaker condition of having H\"{o}lderian growth with $ d = 2 $ (c.f., \cite{fercoq2019adaptive}). Each of these schemes, however, is designed for a narrow family of first-order methods, and typically relies on learning appropriately-accurate approximations of the parameter values characterizing functions in a particular class (e.g., learning the Lipschitz constant $ L $  when $ f $ is assumed to be smooth). There is not a simple, general and easily-implementable meta-heuristic  underlying any of the schemes.
 
Going beyond the setting of $ f $ being smooth and $ d = 2 $, and going beyond being designed to apply to a narrow family of first-order methods, the foremost contributions on restart schemes are due to Roulet and d'Aspremont \cite{roulet2017sharpness}, who consider all $ f $ possessing H\"{o}lderian growth (which they call ``sharpness''), and having H\"{o}lder continuous gradient with exponent $ 0 < \nu \leq 1 $. Their work aims, in part, to make algorithmically-concrete the earlier abstract analysis of Nemirovski and Nesterov \cite{NemNes85}. 

The restart schemes of Roulet and d'Aspremont result in optimal time bounds when particular algorithms are employed in the schemes, assuming scheme  parameters are set to appropriate values -- values that generally would be unknown in practice.  However, for smooth $ f $ (i.e., $ \nu = 1 $), they develop (\cite[\S3.2]{roulet2017sharpness}) an adaptive grid-search procedure within the scheme, to accurately approximate the required values, leading to overall time bounds that differ from the optimal bounds only by a logarithmic factor. Moreover, for general $ f $ for which the optimal value $ f^* $ is known, they show (\cite[\S4]{roulet2017sharpness}) that when an especially-simple restart scheme employs Nesterov's universal fast gradient method \cite{nesterov2015universal}, nearly optimal time bounds result.

The restart schemes with established time bounds either rely on values that would not be known in practice (e.g., $ f^* $), or assume the function $ f $ is of a particular form (e.g., smooth and having H\"{o}lderian growth), and adaptively approximate values characterizing the function (e.g., $ \mu $ and $ d $) in order to set scheme parameters to appropriate values.  By contrast, the restart scheme we propose depends only on how much progress has been made, how much the objective value has been decreased. 

Nonetheless, when Nesterov's universal fast gradient method \cite{nesterov2015universal} is employed in our scheme, it achieves nearly optimal complexity bounds for $ f $ having H\"{o}lderian growth and  H\"{o}lder continuous gradient with exponent $ 0 <   \nu \leq  1 $ (Corollary~\ref{cor.hb}).  Here, our results break new ground in their combination of (1) generality of the class of functions $ f $, (2) attainment of nearly-optimal complexity bounds, and (3) avoidance of relying on problem information that would not be available in practice.

\subsection{Positioning within the literature on nonsmooth convex optimization.} \label{sect.af}

Employing the scheme even with the simple algorithm $ \subgrad $  leads to new and consequential results when $ f $ is a nonsmooth convex function that is Lipschitz  on an open neighborhood of $ Q $. In this regard, we note that for such functions which in addition have growth degree $ d = 1 $, there was considerable interest in obtaining linear convergence even in the early years of research on subgradient methods
 (see \cite{polyak1963gradient,polyak1977subgradient,goffin1977convergence,shor1985minimization} for discussion and references). The goal was accomplished, however, only under substantial assumptions.
 
 Interest in the setting continues today. In recent years, various algorithms have been developed with complexity bounds depending only logarithmically on $ \epsilon $  (c.f., \cite{bolte2017error,gilpin2010first,johnstone2020faster,renegar2016efficient,yang2018rsg}). 
Nevertheless, each of the algorithms for which a logarithmic bound has been established depends on exact knowledge of values characterizing an optimization problem (e.g., $ f^* $), or on accurate estimates of values (e.g., the distance of the initial iterate from optimality, or a ``growth constant''), that generally would be unknown in practice. None of the algorithms entirely avoids the need for special information, although a few of the algorithms are capable of adaptively learning appropriately-accurate approximations to nearly all of the values they rely upon. (Most notable for us, in part due to its generality, is Johnstone and Moulin \cite{johnstone2020faster}.)  

By contrast, given that the algorithm $ \subgrad $ does not rely on parameter values characterizing problem structure, the same is true for our synchronous scheme when $ \subgrad $ is the method employed. It is thus consequential that the resulting time bound, presented in Corollary~\ref{cor.eb}, is proportional to $ \log (1/\epsilon) $. (The total number of subgradient evaluations is proportional to $ \log(1/\epsilon)^2 $.) As mentioned earlier, for general convex optimization problems with Lipschitz objective functions, this provides the first algorithm which both relies on no special information about the optimization problem being solved, and for which it known that if the objective function $ f $ happens to possess linear growth (i.e., happens to be H\"{o}lderian with $ d = 1 $), an $ \epsilon $-optimal solution is computed with a total number of subgradient evaluations depending only logarithmically on $ 1/\epsilon $. (Moreover, as we indicated, when $ \smooth $ is employed in our synchronous scheme, not only is logarithmic dependence on $ 1/\epsilon $ obtained for important convex optimization problems, but also linear dependence on the condition number $ M/\mu $, an even stronger result at least in theory, than when $ \subgrad $ is employed in the scheme.)

Perhaps the most surprising aspect of the paper is that the numerous advances are accomplished with such a simple restart scheme. We now begin detailing the scheme and the elementary theory.

\section{{\bf Assumptions}}  \label{sect.b}

Recall that we consider optimization problems of the form
\begin{equation}  \label{eqn.ba} 
\begin{array}{rl}
 \min & f(x) \\
\textrm{s.t.} & x \in Q \; , 
\end{array} \end{equation}
where $ f $ is a convex function, and $ Q $ is a closed convex set contained in the (effective) domain of $ f $. The set of optimal solutions, $ X^* $, is assumed to be nonempty. Recall $ f^* $ denotes the optimal value.

Recall also that $ \fom $ denotes a first-order method capable of solving some class of convex optimization problems of the form (\ref{eqn.ba}), in the sense that for any problem in the class, when given an initial point $ x_0 $ and accuracy $ \epsilon  > 0 $, the method can be specialized to provide an algorithm to generate a feasible iterate $ x_k $ guaranteed to be an $ \epsilon $-optimal solution. Let $ \fom(\epsilon) $ denote the specialized version of $ \fom $. 

For $ \fom = \subgrad $, we let $ \subgrad( \epsilon) $ denote the first-order method with iterates $ x_{k+1} = P_Q( x_k - \frac{\epsilon}{\| g_k \|^2} g_k ) $, where $ g_k \in \partial f(x_k) $. For $ \fom = \accel $, we let $ \accel( \epsilon) = \accel $, independent of $ \epsilon $. 

\subsection{Assumptions on algorithms.} \label{sect.ba}
In considering a first order method $ \fom $ which requires approximately the same amount of work during each iteration, we assume that the number of iterations (number of oracle calls) sufficient to achieve $ \epsilon $-optimality can be expressed as a function of two specific parameter values, as typically is the case. In particular, we assume there is a function $ K_{\fom}: \mathbb{R}_+ \times \mathbb{R}_+ \rightarrow \mathbb{R}_+ $ satisfying
\begin{equation} \label{eqn.bb} 
  \dist(x_0, X^*) \leq \delta \quad \Rightarrow \quad  \min \{ f(x_k) \mid 0 \leq k \leq K_{\fom}(\delta , \epsilon) \} \, \leq \, f^* + \epsilon.
   \end{equation}  
We also assume the function $ K_{\fom} $ satisfies the following natural condition:
\begin{equation}  \label{eqn.bc}
    \delta \leq \delta' \textrm{ and } \epsilon \geq \epsilon' \quad \Rightarrow \quad K_{\fom}(\delta, \epsilon ) \leq K_{\fom}(\delta', \epsilon') \; . 
    \end{equation}
   
For example, consider the class of problems of the form (\ref{eqn.ba})   for which all subgradients of $ f $ on $ Q $ satisfy $ \| g \| \leq M $ for fixed $ M $ (we slightly abuse terminology and say that $ f $ is ``$ M $-Lipschitz on $ Q $''). An appropriate algorithm is $ \subgrad $, for which it is well known that the function $ K_{\fom} = K_{\subgrad} $ given by 
\begin{equation}   \label{eqn.bd} 
K_{\mathtt{subgrad}}(\delta, \epsilon) :=  (M  \delta/\epsilon )^2  
\end{equation} 
satisfies (\ref{eqn.bb}). (A simple proof is included in Appendix~\ref{sect.j}.) 

As another example, for problems in which $ f $ is $ L $-smooth on an open neighborhood of $ Q $, the first-order method $ \accel $ can be applied, for which  (\ref{eqn.bb})   holds with $ K_{\fom} = K_{\mathtt{accel}} $ given by 
\begin{equation}  \label{eqn.be}  
  K_{\mathtt{accel}}(\delta,\epsilon ) :=   \delta \sqrt{2L/\epsilon}   \; . 
\end{equation}
  
The same function bounds the number of oracle calls for $ \fista $, in the more general setting of proximal methods. (See \cite[Thm 4.4]{beck2009fast} for $ \fista $ -- for the special case of $ \accel $, choose $ g $ in that theorem to be the indicator function for $ Q $.)

A third general example we rely upon pertains to the notion of ``smoothing'' promoted by Nesterov\cite{nesterov2005smooth}, in which a nonsmooth convex objective function $ f $ is replaced by a smooth convex approximation $ f_{\eta} $, depending on parameter $ \eta $ which controls the error in how well $ f_{\eta} $ approximates $ f $, as well as determines the smoothness of $ f_{\eta} $ (i.e., determines the Lipschitz constant $ L_{\eta} $ for which $ \| \nabla f_{\eta}(x) - f_{\eta}(y) \| \leq L_{\eta} \| x - y \| $).  Nesterov showed that for functions $ f $ of particular structure, there is a smoothing $ f_{\eta} $ such that if $ \eta $ is chosen appropriately, and if an accelerated gradient method is applied with $ f_{\eta} $ in place of $ f $, then an iteration bound growing only like $ O(1/\epsilon) $ results for obtaining an iterate $ x_k $ which is an $ \epsilon $-optimal solution for $ f $, the objective function of interest. Compare this with the worst-case $ O(1/ \epsilon^2) $ bound resulting from applying a subgradient method to $ f $.

Nesterov's notion of smoothing was further developed by Beck and Teboulle in \cite{beck2012smoothing}, where they used a slightly more general definition than the following:
\begin{quote}
An  ``$ (\alpha, \beta) $ smoothing of $ f $'' is a family of functions $ f_{\eta} $ parameterized by $ \eta > 0 $, where for each $ \eta $, the function $ f_{\eta} $ is smooth on a neighborhood of $ Q $ and satisfies
\begin{enumerate}
\item $ \| \grad f_{\eta}(x) - \grad f_{\eta}(y) \| \leq \smfrac{\alpha }{\eta }  \| x - y \| \; , \quad \forall x,y \in Q $.
\item $ f(x) \leq f_{\eta}(x) \leq f(x) + \beta \eta \; , \quad \forall x \in Q $,
\end{enumerate}
\end{quote}
(Unlike \cite{beck2012smoothing}, we restrict attention to the Euclidean norm due to our interest in H\"{o}lderian growth.)

For example, for the piecewise-linear convex function 
\begin{equation} \label{eqn.bf} 
  f(x) = \max \{ a_i^T x - b_i \mid i = 1, \ldots, m \} \quad \textrm{where  }  a_i \in \mathbb{R}^n,  b_i \in \mathbb{R} \; , 
  \end{equation}
a $ (\max_i \| a_i \|^2, \ln (m)) $-smoothing\footnote{This is example 4.5 in \cite{beck2012smoothing}, but there the Euclidean norm is not used, unlike here. To verify correctness of the value $ \alpha = \max_i \| a_i \|^2 $, it suffices to show for each $ x $ that the largest eigenvalue of the positive-definite matrix $ \nabla^2 f_{\eta}(x) $ does not exceed $ \max_i \| a_i \|^2/ \eta $. However, \\
 $ \nabla^2 f_{\eta}(x) = - \smfrac{1}{\eta } \grad f_{\eta}(x)  \grad f_{\eta}(x)^T  +  \smfrac{1}{ \eta \sum_i \exp ((a_i^T x + b_i)/ \eta ) } \sum_i \exp((a_i^T x + b_i)/\eta) \, a_i a_i^T \; , $     
and hence the largest eigenvalue of $ \nabla^2 f(x) $ is no greater than $ 1/ \eta $ times the largest eigenvalue of the matrix \\
 $ \smfrac{1}{\sum_i \exp ((a_i^T x + b_i)/ \eta ) } \sum_i \exp((a_i^T x + b_i)/\eta) \, a_i a_i^T $, which is a convex combination of the matrices $ a_i a_i^T $.  Consequently $ \alpha = \max_i \| a_i \|^2 $ is indeed an appropriate choice.}
  is given by 
\begin{equation} \label{eqn.bg}
    f_{\eta}(x) = \eta \ln \left( \sum_{i=1}^m \exp ( (a_i^T x - b_i)/ \eta) \right)  \; . 
\end{equation}

Similarly, for eigenvalue optimization problems involving the maximum eigenvalue function
\[ \begin{array}{rl} \min & \lambda_{\max}(X) \\
\textrm{s.t.} & {\mathcal A}(X) = b \; , \end{array} \]
(where $ {\mathcal A} $ is a linear map from $ \Sym $ ($ n \times n $ symmetric matrices) to $ \mathbb{R}^{m} $), Nesterov\cite{nesterov2007smoothing}  showed with respect to the trace inner product that a $ (1, \ln(n)) $-smoothing is given by
\begin{equation}  \label{eqn.bh} 
    f_{\eta}(x) = \eta \ln \left( \sum_{j=1}^n \exp ( \lambda_j(x) / \eta) \right)  \; . 
\end{equation}
(This smoothing generalizes to all hyperbolic polynomials, $ X \mapsto \det(X)  $ being a special case (\cite{renegar2017accelerated}).)  

In Appendix \ref{sect.k}, we demonstrate that it is straightforward to develop a first-order method $ \smooth $ (based on $ \accel $ (or more generally, $ \fista $)) that makes use of a smoothing to compute an $ \epsilon $-optimal solution for $ f $, subject to $ x \in Q $, within a number of steps (gradient evaluations) not exceeding
\begin{equation}  \label{eqn.bi} 
K_{\smooth}(\delta, \epsilon) := 3 \delta \sqrt{2 \alpha \beta }/ \epsilon \; . 
\end{equation}

When smoothing is possible, the iteration count for obtaining an $ \epsilon $-optimal solution is typically far better than the count for the subgradient method. For example, noting that the piecewise-linear function (\ref{eqn.bf}) is Lipschitz with constant $ M := \max_i \| a_i \| $, we have the upper bounds
\[   K_{\subgrad}(\delta, \epsilon) = (M \delta/ \epsilon)^2 \quad \textrm{and} \quad K_{\smooth}(\delta, \epsilon) = 3 M \delta \sqrt{2 \ln(m)} / \epsilon \; . \]
Similarly for eigenvalue optimization.

However, the cost per iteration is typically much higher for $ \smooth $ than for $ \subgrad $, because rather than computing a subgradient of $ f $, $ \smooth $ requires a gradient of $ f_{\eta} $ (for specified $ \eta $). In the case of the maximum eigenvalue function, for example, a subgradient at $ X $ is simply $ v v^T $, where $ v $ is any unit-length eigenvector for eigenvalue $ \lambda_{\max}(X) $, whereas computing the gradient $ \grad f_{\eta}(X) $ essentially requires a full eigendecomposition of $ X $.

\subsection{Assumptions for theory.} \label{sect.bb}
Throughout the paper, theorems are stated in terms of the function $ K_{\fom} $ and values
\[ 
 D(\hat{f}) = \sup \{ \dist(x,X^*) \mid  x \in Q  \textrm{ and } f(x) \leq \hat{f} \} 
\] 
(assuming $ \hat{f} \geq f^* $). For the theorems to be meaningful, the values $ D( \hat{f}) $ must be finite.

In corollaries throughout the paper, it is assumed additionally that the convex function $ f $ has H\"{o}lderian growth, in that there exist positive constants $ \mu $  and $ d  $ for which
\[  
   x \in Q \textrm{ and } f(x) \leq f( \mathbf{x_0})  \quad \Rightarrow \quad f(x) - f^* \geq \mu \,  \dist(x,X^*)^d \; ,
\] 
where $ \mathbf{x_0} $ is the initial iterate.  (This assumption can be slightly weakened at the expense of heavier notation.\footnote{In particular, analogous results can be obtained under the weaker assumption that  H\"{o}lderian growth holds on a smaller set $ \{ x \in Q \mid \dist(x,X^*) \leq r \} $ for some $ r > 0 $. The key is that due to convexity of $ f $, for $ \hat{f} \geq f^* $ we would then have 
\[     D( \hat{f})  \leq  \begin{cases}  \big( (\hat{f} - f^*)/\mu  \big)^{1/p} & \textrm{if $ \hat{f} - f^* \leq \mu \, r^p $} \\
(\hat{f} - f^*)/(\mu \, r^{p-1}) & \textrm{if $ \hat{f} - f^* \geq \mu \, r^p $} \; .        
\end{cases}    \]}) It is easy to prove that convexity of $ f $ forces $ d $ to satisfy $ d \geq 1 $. It also is not difficult to show that for smooth convex functions possessing H\"{o}lderian growth, it must be the case that $ d \geq 2 $.

For each of our theorems, the corollaries are deduced in straightforward manner from  the immediate fact that if $ f $ has H\"{o}lderian growth and $ \hat{f} \geq f^* $, then 
\[ 
  D(\hat{f}) \leq  \left( \smfrac{\hat{f} - f^*}{\mu} \right)^{1/d} \; \; , \]
and hence for any $ \epsilon > 0 $,    
\begin{equation} \label{eqn.bj} 
  K_{\fom}( D( \hat{f}), \epsilon) \leq K_{\fom} \left( \left( \smfrac{\hat{f} - f^*}{\mu} \right)^{1/d}, \epsilon \right)  
  \end{equation}
  (making use of the natural assumption (\ref{eqn.bc})).
 
\section{{\bf Idealized Setting: When the Optimal Value is Known}} \label{sect.c} 

When the optimal value is known, optimal complexity bounds can readily be established without having to run multiple copies of a first order method in parallel, as we show in the current section. This section is similar in spirit  to the analysis in \cite[\S4]{roulet2017sharpness}, but is done in a manner which immediately applies to numerous first-order methods by focusing on the function $ K_{\fom}(\delta, \epsilon) $. Our development serves to motivate the structure of our restart scheme (presented in \S\ref{sect.d}), and serves to make more transparent the proofs of results regarding our scheme.

The following scheme never terminates, and computes an $ \epsilon $-optimal solution for every $ \epsilon > 0 $. 
\vspace{2mm}

\noindent 
$ \textrm{~} $ \quad {\bf  ``Polyak Restart Scheme''} \\ 
$ \textrm{~} $ \qquad (0) Input: $ f^* $, $ \mathbf{x_0}    \in Q $ (an initial point) \\
$ \textrm{~} $ \qquad \, \, \, \, Initialize: Let $ x_0 = \mathbf{x_0}  $. \\
$ \textrm{~} $  \qquad   (1) Compute  $ \tilde{\epsilon} := \frac{1}{2}( f( x_0) - f^* ) $. \\
$ \textrm{~} $ \qquad    (2) Initiate $ \fom( \tilde{\epsilon}) $ at $ x_0  $, \\
$ \textrm{~} $ \qquad \qquad \quad  and continue until reaching the first iterate $ x_k $ satisfying $ f(x_k) \leq f(x_0)-\tilde{\epsilon} $. \\
$ \textrm{~} $ \qquad    (3) Substitute $ x_0 \leftarrow x_k $, and go to Step 1.  
\vspace{2mm}

We name the scheme after B.T. Polyak due to his early promotion of the subgradient method given by $ x_{k+1} = x_k - \frac{f(x_k) - f^*}{ \| g_k \|^2} g_k $, where the step size $ \alpha_k =\frac{f(x_k) - f^*}{ \| g_k \|^2} $ is updated to be ``ideal'' at every iteration (c.f., \cite{polyak1987introduction}). Our Polyak restart scheme is a bit different, in using the value $ f(x_0) - f^* $ to set a goal that when achieved, results in a restart.  This restart scheme is effective even for first-order methods that do not rely on $ \epsilon $ as input, such as $ \accel $.

The update rule based on $ \tilde{\epsilon} = \frac{1}{2}( f( x_0) - f^* ) $ is slightly arbitrary, in that for any constant $ 0 < \kappa < 1 $, the rule with $ \tilde{\epsilon} = \kappa \, ( f( x_0) - f^* ) $ would lead to iteration bounds that differ only by a constant factor from the bounds we deduce below.

\begin{thm} \label{thm.ca} 
 For each $ \epsilon $ satisfying $ 0 < \epsilon < f( \mathbf{x_0}) - f^*  $, the Polyak restart scheme computes an $ \epsilon $-optimal solution within a total number of iterations (oracle calls) for $ \fom $ not exceeding  
\begin{align}  \label{eqn.ca} 
   \sum_{n=-1}^{\bar{N}}  K_{\fom}(D_n, 2^n \epsilon ) 
\quad \textrm{with }   D_n :=   \min \{ D(f^* + 4 \cdot 2^n \epsilon), D( f( \mathbf{x_0})) \} \; ,  
\end{align}
where $ \bar{N} = \lfloor \log_2 \left(  \frac{f( \mathbf{x_0}) - f^*}{\epsilon} \right) \rfloor - 1 $.          
\end{thm}
\noindent {\bf Proof:} Denote the sequence of values $ \tilde{\epsilon} $ computed by the scheme as $ \tilde{\epsilon}_1, \tilde{\epsilon}_2, \ldots $. A simple inductive argument shows the sequence is infinite, by observing that when working with $ \tilde{\epsilon}_i $, the scheme employs $ \fom( \tilde{\epsilon}_i ) $, which is guaranteed to compute an $ \tilde{\epsilon}_i $-optimal solution, at which time $ \tilde{\epsilon}_i $ is replaced by a value $ \tilde{\epsilon}_{i+1} \leq \tilde{\epsilon}_i/2 $.

For each $ i $, let $ y_{i,0}, \ldots, y_{i,k_i} $ denote the sequence of iterates for $ \fom( \tilde{\epsilon}_i) $ -- thus, $ \tilde{\epsilon}_i = \frac{1}{2}(f(y_{i,0}) - f^*) $ and $ y_{i,k_i} $ is the first of these iterates which is an $ \tilde{\epsilon}_i $-optimal solution. Consequently, 
\begin{equation}  \label{eqn.cb} 
      k_i \leq K_{\fom} (\, D(f(y_{i,0})), \tilde{\epsilon}_i \, ) \, \, = K_{\fom} (\, D( f^* + 2 \tilde{\epsilon}_i), \tilde{\epsilon}_i \, ) \; . \end{equation}

Fix any $ \epsilon $ satisfying $ 0 < \epsilon < f( \mathbf{x_0}) - f^*  $. For each $ i $, let $ n_i $ be the largest integer satisfying $ 2^{n_i} \epsilon \leq \tilde{\epsilon}_i $. Of course we also have $ 4 \cdot 2^{n_i} \epsilon > 2 \tilde{\epsilon}_i $. Consequently, by (\ref{eqn.cb}), 
\begin{equation}  \label{eqn.cc} 
 k_i \leq K_{\fom}( \, D(f^* + 4 \cdot 2^{n_i} \epsilon), 2^{n_i} \epsilon \, ) \; . 
\end{equation}

Due to the relations $ \tilde{\epsilon}_{i+1} \leq \frac{1}{2} \tilde{\epsilon}_i $, the sequence $ n_1, n_2, \ldots $ is strictly decreasing. Moreover, $ n_1 $ is precisely the integer $ \bar{N} $ in the statement of the theorem.   

Let $ i' $ be the last index $ i $ for which $ \tilde{\epsilon}_i > \epsilon/2 $. Then $ f(y_{i'+1,0}) - f^* = 2 \tilde{\epsilon}_{i' + 1} \leq \epsilon $ -- thus, since $ y_{i'+1,0} = y_{i',k_{i'}} $, we have that $y_{i',k_{i'}} $ is an $ \epsilon $-optimal solution. Moreover, since $ \tilde{\epsilon}_{i'} > \epsilon/2 $, we have $ n_{i'} \geq -1 $.

In all, the total number of iterations of $ \fom $ required by the Polyak scheme to compute an $ \epsilon $-optimal solution does not exceed
\[    \sum_{i=1}^{i'} k_i \leq \sum_{n=-1}^{\bar{N}} K_{\fom} \big(   D(f^* + 4 \cdot 2^n \epsilon), 2^n \epsilon \big)  \; , \]
establishing the theorem. \hfill $ \Box $

\begin{cor} \label{cor.cb} 
Assume $ f $ has H\"{o}lderian growth.
 For each $ \epsilon $ satisfying $ 0 < \epsilon  < f( \mathbf{x_0}) - f^*  $, the Polyak scheme computes an $ \epsilon $-optimal solution within a total number of iterations of $ \fom $ not exceeding  
\begin{equation}  \label{eqn.cd}
  \sum_{n=-1}^{ \bar{N}} K_{\fom} \left( \left( \smfrac{4 \cdot 2^n \epsilon }{\mu} \right)^{1/d}, 2^n \epsilon \right) \; , 
  \end{equation} 
where $ \bar{N} = \lfloor \log_2 \left(  \frac{f( \mathbf{x_0}) - f^*}{\epsilon} \right) \rfloor - 1 $. 
\end{cor}
\noindent {\bf Proof:}  Simply substitute into Theorem~\ref{thm.ca}  using (\ref{eqn.bj}). \hfill $ \Box $

\begin{cor}  \label{cor.cc} 
Assume $ f $ is $ M $-Lipschtiz on $ Q $, and has H\"{o}lderian growth. For $ \epsilon $ satisfying $ 0 < \epsilon < f( \mathbf{x_0}) - f^* $, the Polyak scheme with $ \fom = \subgrad $ computes an $ \epsilon $-optimal solution within a total number of $ \subgrad $ iterations $ K $, where
\begin{align}
  d & = 1 \, \,   \Rightarrow \, \, K \leq \,  (\bar{N}+2) (4 M/ \mu)^2 \; ,  \label{eqn.ce} \\
  d &> 1 \, \,  \Rightarrow \, \,  K \leq \,  \left(  \frac{ 4^{1/d} M }{\mu^{1/d} \, \epsilon^{1 - 1/d}} \right)^2 \min \left\{ \frac{ 16^{1 - 1/d}}{4^{1 - 1/d} - 1}, \, \bar{N} + 5 \right\} \; , \label{eqn.cf} 
  \end{align}
with $ \bar{N} = \lfloor \log_2 \left(  \frac{f( \mathbf{x_0}) - f^*}{\epsilon} \right) \rfloor - 1 $.  
 \end{cor}
\noindent {\bf Proof:}  
From $ K_{\subgrad}(\delta, \bar{\epsilon}  ) :=   (M \delta/ \bar{\epsilon} )^2  $ follows
\[  
    K_{\subgrad} \left( (4 \cdot 2^n \epsilon/\mu)^{1/d}, \, 2^n \epsilon \right)  \leq  \left(  \frac{M (4 \cdot 2^n \epsilon/ \mu )^{1/d}}{2^n \epsilon} \right)^2 = \left( \frac{ 4^{1/d} M }{\mu^{1/d} \epsilon^{1 - 1/d}} \right)^2   \left( \frac{1}{ 4^{1 - 1/d}} \right)^n \; , 
\]
and hence (\ref{eqn.cd})   is bounded above by
\[ 
 \left( \frac{ 4^{1/d} M }{\mu^{1/d} \epsilon^{1 - 1/d}} \right)^2  \sum_{n=-1}^{\bar{N}} \left( \frac{1}{ 4^{1 - 1/d}} \right)^n \; . 
\]  The implication for $ d = 1 $ is immediate by Corollary~\ref{cor.cb}. On the other hand, since for $ d > 1 $, 
\[  
  \sum_{n=-1}^{\bar{N}} \left( \frac{1}{ 4^{1 - 1/d}} \right)^n < 
 \min \left\{ \sum_{n=-1}^{\infty} \left( \frac{1}{ 4^{1 - 1/d}} \right)^n, \, \, \bar{N} + 5  \right\} 
 =   \min 
  \left\{ \frac{ 16^{1 - 1/d}}{4^{1 - 1/d} - 1}, \, \bar{N} + 5 \right\} \; ,  
\]  
the implication (\ref{eqn.cf})  is immediate, too. 
 \hfill $ \Box $
\vspace{2mm}

\noindent {\bf Remark:} Note that when $ d = 1 $, as it is for the convex piecewise-linear function (\ref{eqn.bf}), the corollary shows the Polyak scheme with $ \fom = \subgrad $ converges linearly, i.e., the number of iterations grows only like $ \log(1/\epsilon) $ as $ \epsilon \rightarrow 0 $.

The next corollary regards $ \accel $, the accelerated method, for which $ f $ is assumed to be smooth, and thus $ d $, the degree of growth, satisfies $ d \geq 2 $. If the term ``iterations'' is replaced by ``oracle calls,'' the corollary holds for $ \fista $, of which $ \accel $ is a special case, as we mentioned in \S\ref{sect.aa}.

\begin{cor}  \label{cor.cd} 
Assume $ f $ is $ L $-smooth on a neighborhood of $ Q $, and has H\"{o}lderian growth. For $ \epsilon $ satisfying $ 0 < \epsilon < f( \mathbf{x_0}) - f^* $, the Polyak scheme with $ \fom = \accel $ computes an $ \epsilon $-optimal solution within a total number of $ \accel $ iterations $ K $, where
\begin{align}
  d & = 2  \, \,   \Rightarrow \, \,   K \leq  2 (\bar{N}+2) \sqrt{ 2 L/\mu}  \; ,     \label{eqn.cg} \\
  d &> 2 \, \,  \Rightarrow \, \,  K \leq  \frac{ (4/ \mu)^{1/d} \sqrt{2L}}{\epsilon^{ \frac{1}{2} - \frac{1}{d}}}   \min \left\{ \frac{4^{ \frac{1}{2} - \frac{1}{d}}}{2^{\frac{1}{2} - \frac{1}{d}} - 1}, \,  \bar{N} + 3 \right\}  \; , \label{eqn.ch} 
  \end{align}
with $ \bar{N} = \lfloor \log_2 \left(  \frac{f( \mathbf{x_0}) - f^*}{\epsilon} \right) \rfloor - 1 $.    
    \end{cor}
\noindent {\bf Proof:}  
From $ K_{\accel}(\delta, \bar{\epsilon}   ) :=    \delta \sqrt{2 L/ \bar{\epsilon}  }   $ follows
\[  
     K_{\accel}( (4 \cdot 2^n \epsilon/ \mu)^{1/d}, 2^n \epsilon)  \leq  \frac{ \sqrt{2L}  (4 \cdot 2^n \epsilon/ \mu)^{1/d}}{\sqrt{2^n \epsilon }}  = \frac{ \sqrt{2L} (4/ \mu)^{1/d}}{ \epsilon^{ \frac{1}{2} - \frac{1}{d}}} \left( \frac{1}{2^{\frac{1}{2} - \frac{1}{d}}} \right)^n  \; ,              
\] 
and hence (\ref{eqn.cd}) is bounded above by 
\[  
  \frac{ \sqrt{2L} (4/ \mu)^{1/d}}{ \epsilon^{ \frac{1}{2} - \frac{1}{d}}}  \sum_{n = -1}^{\bar{N}}  \left( \frac{1}{2^{\frac{1}{2} - \frac{1}{d}}} \right)^n \; . 
\] 
The implication for $ d = 2 $  is immediate from Corollary~\ref{cor.cb}. On the other hand, since for $ d > 2 $, 
\[    
   \sum_{n=-1}^{\bar{N}} \left(  \frac{1}{ 2^{\frac{1}{2} - \frac{1}{d}}  } \right)^n < \min \left\{ \sum_{n=-1}^{\infty } \left(  \frac{1}{ 2^{\frac{1}{2} - \frac{1}{d}}}\right)^n, \, \bar{N} + 3 \right\} 
            = \min \left\{ \frac{4^{ \frac{1}{2} - \frac{1}{d}}}{2^{\frac{1}{2} - \frac{1}{d}} - 1}, \,  \bar{N} + 3 \right\} \; ,
\] 
the implication (\ref{eqn.ch})  is immediate, too.  \hfill $ \Box $

\noindent {\bf Remark:} Note that when $ d = 2 $ -- a class of functions in which the strongly convex functions form a small subset -- the corollary shows the Polyak scheme with $ \fom = \accel $ converges linearly.

The final corollary of this section regards the Polyak scheme with $ \fom = \smooth $, appropriate when $ f $ is a nonsmooth convex function for which a smoothing is known.

\begin{cor}  \label{cor.ce} 
Assume $ f $ has an $ (\alpha,\beta) $-smoothing, and assume $ f $ has H\"{o}lderian growth.  For $ \epsilon $ satisfying $ 0 < \epsilon < f( \mathbf{x_0}) - f^* $, the Polyak scheme with $ \fom = \smooth $ computes an $ \epsilon $-optimal solution for $ f $ within a total number of $ \smooth $ iterations $ K $, where
\begin{align}
  d & = 1  \, \,   \Rightarrow \, \,   K \leq   12 (\bar{N}+2) \sqrt{ 2 \alpha \beta } / \mu  \; ,     \label{eqn.ci} \\
  d & > 1 \, \,  \Rightarrow \, \,  K \leq  \frac{ 3  \sqrt{2 \alpha \beta} (4/ \mu)^{1/d}}{\epsilon^{ 1 - 1/d}}   \min \left\{ \frac{4^{ 1 - 1/d}}{2^{1 - 1/d} - 1}, \,  \bar{N} + 3 \right\}  \; , \label{eqn.cj} 
  \end{align}
with $ \bar{N} = \lfloor \log_2 \left(  \frac{f( \mathbf{x_0}) - f^*}{\epsilon} \right) \rfloor - 1 $.    
    \end{cor}
\noindent {\bf Proof:}  
From (\ref{eqn.bi})  follows
\[  
		 K_{\smooth}( (4 \cdot 2^n \epsilon/ \mu)^{1/d}, 2^n \epsilon)  \leq \frac{3 \sqrt{2 \alpha \beta} (4 \cdot 2^n \epsilon/\mu)^{1/d}}{2^n \epsilon}  = \frac{3 \sqrt{2 \alpha \beta } (4/\mu)^{1/d}}{\epsilon^{1 - 1/d}} \left( \frac{1}{2^{1 - 1/d } }\right)^n             \; ,              
\] 
and hence (\ref{eqn.cd})  is bounded above by 
\[  
  \frac{ 3 \sqrt{2 \alpha \beta } (4/ \mu)^{1/d}}{ \epsilon^{ 1 - 1/d}}  \sum_{n = -1}^{\bar{N}}  \left( \frac{1}{2^{1 - 1/d}} \right)^n \; . 
\] 
The implication for $ d = 1 $  is immediate from Corollary~\ref{cor.cb}. On the other hand, since for $ d > 1 $, 
\[    
   \sum_{n=-1}^{\bar{N}} \left(  \frac{1}{ 2^{1 - 1/d  }} \right)^n < \min \left\{ \sum_{n=-1}^{\infty } \left(  \frac{1}{ 2^{1 - 1/d}}\right)^n, \, \bar{N} + 3 \right\} 
            = \min \left\{ \frac{4^{ 1 - 1/d}}{2^{1 - 1/d} - 1}, \,  \bar{N} + 3 \right\} \; ,
\] 
the implication (\ref{eqn.cj})  is immediate, too.  \hfill $ \Box $
 
To understand the relevance of the corollary, consider the piecewise-linear convex function $ f $ defined by (\ref{eqn.bf}), which has $ ( \max_i \| a_i \|^2, \ln (m)) $-smoothing given by (\ref{eqn.bg}). The function $ f $ has linear growth (i.e., $ d = 1 $), for some $ \mu > 0 $. Consequently, according to Corollary~\ref{cor.ce}, the Polyak scheme with $ \fom = \smooth $ obtains an $ \epsilon $-optimal solution within
\[  12 ( \bar{N} + 2) \sqrt{2 \ln(m)} \max_i \| a_i \|  / \mu  \textrm{ iterations, where } \bar{N} = \lfloor \log_2 \left(  \smfrac{f( \mathbf{x_0}) - f^*}{\epsilon} \right) \rfloor - 1 \; . \]
On the other hand, $ f $  is clearly Lipschitz with constant $ M = \max_i \| a_i \| $, and thus from Corollary~\ref{cor.cc}, the Polyak scheme with $ \fom = \subgrad $ obtains an $ \epsilon $-optimal solution within
\[ 16 (\bar{N}+2) (\max_i \| a_i \|/ \mu)^2 \; \textrm{ iterations}  . \]

For convex functions which are Lipschitz and have linear growth, the ratio $ M/\mu $ is sometimes referred to as the ``condition number'' of $ f $, and is easily shown to have value at least $ 1 $. While the Polyak scheme with either $ \smooth $ or $ \subgrad $ results in linear convergence, the bound for $ \subgrad $ has a factor depending quadratically on the condition number, whereas the dependence for $ \smooth $ is linear. 

Similarly, for an eigenvalue optimization problem, the dependence on $ \epsilon $ for the Polyak scheme with $ \smooth $ is always at least as good as for the scheme with $ \subgrad $ (and is better  when $ d > 1 $), and for $ \smooth $ the dependence on the condition number $ 1/\mu $ is linear whereas for $ \subgrad $ it is quadratic.

Nonetheless, while applying $ \smooth $ typically results in improved iteration bounds, the cost per iteration of $ \smooth $ can be considerably higher than the cost per iteration of $ \subgrad $, as was discussed at the end of \S\ref{sect.ba}  for the eigenvalue optimization problem.

In the following sections where $ f^* $ is not assumed to be known, the computational scheme we develop has nearly the same characteristics as the Polyak scheme. In particular, the analogous iteration bounds for the various choices of $ \fom $ are only worse by a logarithmic factor, to account for multiple instances of $ \fom $ making iterations simultaneously.

\section{{\bf Synchronous Parallel Scheme ($ \sparfom $)}} \label{sect.d} 
Henceforth, we assume $ f^* $ is unknown, and hence the Polyak Scheme is unavailable due to its reliance on the value $ \tilde{\epsilon } :=  \frac{1}{2}( f(x_0) - f^*) $. Our approach relies on running several copies of a first-order method in parallel, that is, running algorithms $ \fom( \epsilon_n) $ in parallel, with various values for $ \epsilon_n $. Similar to the Polyak Scheme, each of these parallel algorithms is searching for an iterate $x_k$ satisfying $f(x_k) \leq f(x_0) - \epsilon_n $.

We assume the user specifies an integer $ N \geq 0 $, with $ N + 2 $ being the number of copies of $ \fom $ that will be run in parallel. We also assume the user specifies a value $ \epsilon > 0 $, the goal being to obtain an $ \epsilon $-optimal solution. Ideally the user would choose $ N \approx \log_2 ( \frac{f( \mathbf{x_0}) - f^*}{ \epsilon} )$ where $  \mathbf{x_0} $ is a feasible point at which the algorithms are initiated. But since $ f^* $ is not assumed to be known, our analysis must apply to various choices of $ N $. Our results yield interesting iteration bounds even for the easily computable choice $ N = \max \{ 0, \lceil   \log_2(1/\epsilon) \rceil \} $.   

The algorithms to be run in parallel are $ \fom( \epsilon_n) $ with $ \epsilon_n = 2^n \epsilon $ ($ n = -1, 0, \ldots, N $).\footnote{The choice of values $ 2^n \epsilon $ ($ n = -1, 0, \ldots, N $) is slightly arbitrary. Indeed, given any scalar $ \gamma > 1 $, everywhere replacing $ 2^n \epsilon $ by $ \gamma^n \epsilon $ leads to similar theoretical results.  However, our analysis shows relying on powers of $ 2 $ suffices to give nearly optimal guarantees, although a more refined analysis, in the spirit of \cite{roulet2017sharpness}, might reveal desirability of relying on values of $ \gamma $ other than 2, depending on the setting.}
Notice that if $N>\log_2 ( \frac{f( x_0) - f^*}{ \epsilon} )$, one of these parallel algorithms must have $ \epsilon_n$ within a factor of $2$ of the idealized choice used by the Polyak Restarting Scheme of $\tilde{\epsilon} := \frac{1}{2}(f(x_0)-f^*)$. Using the key ideas posited in \S\ref{sect.ac}, we carefully manage when these parallel algorithms restart, thereby matching the performance of the Polyak Restarting Scheme (up to a small constant factor), even without knowing $f^*$. The remaining details for specifying this restarting scheme given below are secondary to the key ideas from \S\ref{sect.ac}. It is possible to specify details differently and still obtain similar theoretical results.

To ease notation, we write $ \fom_n $ rather than $ \fom(2^n \epsilon) $ or  $ \fom(\epsilon_n) $. Thus, $ \fom_n $ is a version of $ \fom $ designed to compute a $ 2^n \epsilon $-optimal solution, and it does so within $ K_{\fom}(\delta, 2^n \epsilon) $ iterations, where $ \delta = \dist(x_0, X^*) $, with $ x_0 $ being the initial point for $ \fom_n $. 

In this section we specify details for a scheme in which all of the algorithms $ \fom_n $ make an iteration simultaneously. We dub this the ``synchronous parallel first-order method,'' and denote it as $ \sparfom $. (We will note the simple manner in which $ \sparfom $ can also be implemented sequentially, with the algorithms $ \fom_n $ taking turns in making an iteration, in the cyclic order $ \fom_N, \fom_{N-1}, \ldots, \fom_{-1}, \fom_N, \fom_{N-1}, \,  \ldots $.)

\subsection{Details of $ \sparfom $} \label{sect.da} 

We assume each $ \fom_n $ has its own oracle, in the sense that given $ x $, $ \fom_n $ is capable of computing -- or can readily access -- the information it needs for performing an iteration at $ x $, without having to wait in a queue with other copies $ \fom_m $ in need of information. 

We speak of ``time periods,'' each having the same length, one unit of time. The primary effort for $ \fom_n $ in a time period is to make one iteration\footnote{The synchronous parallel scheme is inappropriate for first order methods that have wide variation in the amount of work done in iterations. Instead, the asynchronous parallel scheme (\S\ref{sect.g}) is appropriate.}, the amount of ``work'' required being measured as the number of calls to the oracle.

At the outset, each of the algorithms $ \fom_n $ ($ n = -1, 0, \ldots, N $) is initiated at the same feasible point, $ \mathbf{x_0} \in Q $. 

The algorithm $ \fom_{N} $ differs from the others in that it is never restarted. It also differs in that it has no ``inbox'' for receiving messages. For $ n < N $, there is an inbox in which $ \fom_n $ can receive messages from $ \fom_{n+1} $. 

At any time, each algorithm has a ``task.'' Restarts occur only when tasks are accomplished. When a task is accomplished, the task is updated.

For all $ n $, the initial task of $ \fom_n $ is to obtain a point $ x $ satisfying $ f(x) \leq f( \mathbf{x_0}) - 2^n \epsilon $. Generally for $ n < N  $, at any time, $ \fom_n $ is pursuing the task of obtaining  $ x $ satisfying $ f(x) \leq f(x_{n,0}) - 2^n \epsilon $, where $ x_{n,0} $ is the most recent (re)start point for $ \fom_n $. Likewise for $ \fom_N $ (the algorithm that never restarts), except that $ x_{n,0} $ is replaced by $ \bar{x}_N $, the most recent ``designated point'' in $ \fom_N $'s sequence of iterates. (The first designated point is $ \bar{x}_N = \mathbf{x_0} $.) 

In a time period, the following steps are made by $ \fom_{N} $:  If the current iterate fails to satisfy the inequality in the task for $ \fom_{N} $, then $ \fom_{N} $ makes one iteration, and does nothing else in the time period. On the other hand, if the current iterate, say $ x_{N,k} $,  does satisfy the inequality, then (1) $ x_{N,k} $ becomes the new designated point ($ \bar{x}_N \leftarrow x_{N,k} $) and $ \fom_N $'s task is updated accordingly, (2) the new designated point is sent to the inbox of $ \fom_{N-1} $ to be available  for $ \fom_{N-1} $ at the beginning of the next time period, and (3) $ \fom_N $ computes the next iterate $ x_{N,k+1} $  (without having restarted).

The steps made by $ \fom_n $ ($ n < N $) are as follows:
First the objective value for the current iterate of $ \fom_n $ is examined, and so is the objective value of the point in the inbox (if the inbox is not empty). If the smallest of these function values -- either one or two values, depending on whether the inbox is empty -- does not satisfy the inequality in the task for $ \fom_n $, then $ \fom_n $ completes its effort in the time period by simply clearing its inbox (if it is not already empty) and making one iteration, without restarting.  On the other hand, if the smallest of the function values does satisfy the inequality, then: (1) $ \fom_n $ is restarted at the point with the smallest function value -- this point becomes the new $ x_{n,0} $ in the updated task for $ \fom_n $ of obtaining $ x $ satisfying $ f(x) \leq f(x_{n,0}) - 2^n \epsilon $ -- and makes one iteration, (2) the new $ x_{n,0} $ is sent to the inbox of $ \fom_{n-1} $ (assuming $ n > -1 $) to be available for $ \fom_{n-1} $ at the beginning of the next time period, and (3) the inbox of $ \fom_n $ is cleared.

This concludes the description of $ \sparfom $.
\vspace{2mm}

\subsection{Remarks} \label{sect.db} 
\begin{itemize}
\item Ideally for each $ n > -1 $,  there is apparatus dedicated solely to the sending of messages from $ \fom_n $ to $ \fom_{n-1} $. If messages from all $ \fom_n $ are handled by a single server, then the length of periods might be increased to more than one unit of time, to reflect possible congestion.
\item In light of $ \fom_n $ sending messages only to $ \fom_{n-1} $, the scheme can naturally be made sequential, with $ \fom_n $ performing an iteration, followed by $ \fom_{n-1} $ (and with $ \fom_N $ following $ \fom_{-1} $). 
\item If instead of $ \fom_n $ sending messages only to $ \fom_{n-1} $, it is allowed that after each iteration, the best iterate among all of the copies $ \fom_n $ is sent to the mailboxes of all of the copies, the empirical results of applying $ \sparfom $ can be noticeably improved (see \S~\ref{sect.i}), although our proof techniques for worst-case iteration bounds can then be improved by only a constant factor. Consequently, in developing theory, we consider only the restrictive situation in which $ \fom_n $ sends messages only to $ \fom_{n-1} $, i.e., our theory considers only the situation in which communication is minimized. Even so, our simple theory breaks new ground for iteration bounds, in multiple regards, as will be noted.
\end{itemize}

\section{{\bf  Theory for $ \sparfom $}}   \label{sect.e}

Here we state the main theorem regarding $ \sparfom $, and present corollaries for $ \subgrad $, $ \accel $ and $ \smooth $. The theorem is proven in  \S\ref{sect.f}.  

Our theorem states bounds on the complexity of $ \sparfom $ depending on whether the positive integer $ N $ is chosen to satisfy a certain inequality involving the difference $ f(\mathbf{x_0}) - f^* $. Since we are not assuming $ f^* $ is known, we cannot assume $ N $ is chosen to satisfy the inequality. 

\begin{thm}  \label{thm.ea}  
If the positive integer $ N $ happens to satisfy $ f(\mathbf{x_0}) - f^* < 5 \cdot 2^N \epsilon $, then $ \sparfom $ computes an $ \epsilon $-optimal solution within time
\begin{align} 
  & \bar{N}+ 1 +  3 \sum_{n=-1}^{\bar{N}}  K_{\fom}\left(   D_n , 2^n \epsilon \right)   \label{eqn.ea}  \\   
& \qquad \qquad   \textrm{with } \,    D_n :=   \min \{ D(f^* + 5 \cdot 2^n \epsilon), D( f( \mathbf{x_0})) \}  \; , \nonumber 
\end{align}   
where $ \bar{N} $ is the smallest integer satisfying both $ f(\mathbf{x_0}) - f^* < 5 \cdot 2^{\bar{N}} \epsilon $ and $ \bar{N} \geq -1 $.

In any case, $ \sparfom $ computes an $ \epsilon $-optimal solution within time
\begin{equation} \label{eqn.eb} 
\mathbf{T}_N + 
  K_{\fom}\big( \dist( \mathbf{x_0}, X^*), 2^{N} \epsilon \big)   \; , 
\end{equation} 
 where $ \mathbf{T}_N  $ is the quantity obtained by substituting $ N $ for $ \bar{N} $ in (\ref{eqn.ea}).
\end{thm}

The theorem is proven in \S\ref{sect.f}.

The theorem gives insight in to how the two user-specified parameters $N$ and $\epsilon$ affect performance. If these choices happen to satisfy $ f(\mathbf{x_0}) - f^* < 5 \cdot 2^N \epsilon $, the time bound~(\ref{eqn.ea}) applies, reflecting the fact (revealed in the proof) that the copies $ \fom_{\bar{N}+1}, \ldots , \fom_{N} $ play no essential role in computing an $ \epsilon $-optimal solution. 

On the other hand, the time bound (\ref{eqn.eb})  always holds, even if the inequality $ f(\mathbf{x_0})-f^*<5\cdot2^N \epsilon $ is not satisfied. By choosing $ N $ to be the easily computable value $ N = \max \{ -1, \lceil \log_2 (1/\epsilon) \rceil \} $, the additive term $ K_{\fom}( \dist(\mathbf{x_0}, X^* ), 2^N \epsilon ) $ becomes bounded by the constant $  K_{\fom}( \dist(\mathbf{x_0}, X^*), 1) $ regardless of the value $ \epsilon $. Thus, with this choice of $ N $, understanding the dependence of the iteration bound on $ \epsilon $ shifts to focus on the value $ \mathbf{T}_N $.     (The point of $ \sparfom $ never restarting $ \fom_N $ is precisely to ensure that the additive term can be made independent of $ \epsilon $ by using an easily-computable choice for $ N $ which grows only like $ \log(1/\epsilon) $ as $ \epsilon \rightarrow 0 $.)

As there are $ N + 2 $ copies of $ \fom $ being run in parallel, the total amount of work (total number of oracle calls) is proportional to the time bound multiplied by $ N +2 $. Of course the same bound on the total amount of work applies if the scheme is performed sequentially rather than in parallel (i.e., $ \fom_n $ performs an iteration, followed by $ \fom_{n-1} $, with $ \fom_N $ following $ \fom_{-1} $). 
\vspace{2mm}

The bound (\ref{eqn.ea})  of Theorem~\ref{thm.ea}  is of the same form as the bound (\ref{eqn.ca})  of Theorem~\ref{thm.ca}.  It is thus no surprise that with proofs exactly similar to the ones for Corollaries \ref{cor.cc}, \ref{cor.cd}  and \ref{cor.ce}, specific time bounds can be obtained for when $ \sparfom $ is implemented with $ \subgrad $, $ \accel $ and $ \smooth $, respectively.

\begin{cor}  \label{cor.eb} 
Assume  $ \fom $ is chosen as $ \subgrad $, $ \accel $ or $ \smooth $, and make the same assumptions as in Corollary \ref{cor.cc}, \ref{cor.cd}  or \ref{cor.ce}, respectively (except do not assume $ f^* $ is known). 

If $ \sparfom $ is applied with $ \fom $, and if $ \epsilon $ and $ N $ happen to satisfy $ f( \mathbf{x_0}) - f^* < 5 \cdot 2^N \epsilon $, then letting $ T $ denote the time sufficient for obtaining an $ \epsilon $-optimal solution, we have

\begin{tabular}{cccl} \\
\multirow{2}{*}{$\subgrad$} 
&  $ d = 1 $  & $ \Rightarrow $ & $ T \leq \, \bar{N}+1  + 3 (\bar{N}+2) (5  M/ \mu)^2 \; , $  \\
 & $ d > 1 $ & $ \Rightarrow $ & $ T \leq \, \bar{N}+1 +  3 \left(  \frac{ 5^{1/d} M }{\mu^{1/d} \, \epsilon^{1 - 1/d}} \right)^2 \min \left\{ \frac{ 16^{1 - 1/d}}{4^{1 - 1/d} - 1}, \, \bar{N} + 5 \right\} \; ,  $ \\  \\
\multirow{2}{*}{$\accel$} 
& $ d = 2 $ & $ \Rightarrow $ & $ T \leq  \bar{N} + 1 + 3 (\bar{N}+2) \sqrt{ 10 L/\mu} \; ,  $ \\
& $ d > 2 $ & $ \Rightarrow $ & $ T \leq \, \bar{N} + 1 + \frac{3 (5/ \mu)^{1/d} \sqrt{2L}}{\epsilon^{ \frac{1}{2} - \frac{1}{d}}}   \min \left\{ \frac{4^{ \frac{1}{2} - \frac{1}{d}}}{2^{\frac{1}{2} - \frac{1}{d}} - 1}, \,  \bar{N} + 3 \right\} \; , $ \\ \\
\multirow{2}{*}{$\smooth$}
& $ d = 1 $ & $ \Rightarrow $ & $ T \leq \bar{N} + 1 +  45 (\bar{N}+2) \sqrt{ 2 \alpha \beta } / \mu \; , $ \\
& $ d > 1 $ &  $ \Rightarrow $ & $ T \leq \bar{N}+ 1 +  \frac{ 9  \sqrt{2 \alpha \beta} (5/ \mu)^{1/d}}{\epsilon^{ 1 - 1/d}}   \min \left\{ \frac{4^{ 1 - 1/d}}{2^{1 - 1/d} - 1}, \,  \bar{N} + 3 \right\} \; .  $ \\ ~
\end{tabular} 

In any case, a time bound is obtained by substituting $ N $ for $ \bar{N} $ above, and adding
\begin{align} 
 ( M \, \dist( \mathbf{x_0}, X^*)/(2^N \epsilon))^2  & \quad  \textrm{for $ \subgrad \; , $} \label{eqn.ec}  \\
 \dist( \mathbf{x_0}, X^*) \sqrt{L/ (2^{N-1} \epsilon)} & \quad \textrm{for $  \accel \; , $} \nonumber \\
  3 \, \dist( \mathbf{x_0}, X^*) \sqrt{2 \alpha \beta}/(2^N \epsilon) & \quad \textrm{for $ \smooth \; .  $} \nonumber
\end{align}
\end{cor}

The bound for $ \subgrad $ when $ d = 1 $ and $ N = \max \{ 0, \lceil \log_2(1/\epsilon) \rceil \} $ establishes $ \sparfom $ as the first algorithm for nonsmooth convex optimization which both has total work proven to depend only logarithmically on $ 1/\epsilon $, and which relies on nothing other than being able to compute subgradients (in particular, does not rely on special knowledge about $ f $ (such as knowing the optimal value, or knowing that $ f $ is H\"{o}lderian with $ d = 1 $, etc.)).  

Similarly, but now for $ d = 2 $,  the easily computable choice $ N = \max \{ 0, \lceil \log_2(1/ \epsilon) \rceil \} $ results for $ \accel $ in a time bound that grows only like $ \sqrt{L/ \mu} \log(1/\epsilon) $ as $ \epsilon \rightarrow 0 $, although the total amount of work grows like $ \sqrt{L/\mu} \log(1/\epsilon)^2 $. This upper bound on the total amount of work (number of gradient evaluations) is within a factor of $ \log(1/\epsilon) $ of the well-known lower bound for strongly-convex smooth functions (\cite{nemirovsky1983problem}), a small subset of the functions having H\"{o}lderian growth with $ d = 2 $.  For $ d > 2 $, the work bound also is within a factor of $ \log(1/ \epsilon) $ of being optimal, according to the lower bounds stated by Nemirovski and Nesterov \cite[page 26]{NemNes85}. 

Finally, $ \sparfom $ with $ \fom = \smooth $ and $ N = \max \{ 0, \lceil \log_2(1/\epsilon) \rceil \} $ shares all of the strengths described at the end of \S\ref{sect.c}  for the Polyak Restart Scheme with $ \fom = \smooth $ (and has the additional property that $ f^* $ need not be known.) Particularly notable is, when $ d = 1 $, the combination of linear dependence on the condition number (for the examples considered in \S\ref{sect.c}, among others), and work bound depending only logarithmically on $ 1/\epsilon $. For no other smoothing algorithm has logarithmic dependence on $ 1/\epsilon $ been established in general when $ d = 1 $ (other than our Polyak scheme with $ \fom = \smooth $). Likewise, when $ d > 1 $, no other smoothing algorithm has been shown to have as good of dependence on $ \epsilon $.

\noindent 
{\bf Proof of Corollary~\ref{cor.eb}:} The proof is exactly similar to the proofs of Corollaries \ref{cor.cc}, \ref{cor.cd}  and \ref{cor.ce}. We focus on $ \subgrad $ to show the minor changes needing to be made to the proof of Corollary~\ref{cor.cc}. 

Just as Theorem~\ref{thm.ca}  immediately yielded the iteration bound (\ref{eqn.cd})  used in the proofs of Corollaries~\ref{cor.cc} , \ref{cor.cd}  and \ref{cor.ce}, so does Theorem~\ref{thm.ea}  yield the time bound 
\begin{equation}  \label{eqn.ed} 
  \bar{N} + 1 + 3 \sum_{n=-1}^{\bar{N}} K_{\fom} \left( \, \left( \smfrac{5 \cdot 2^n \epsilon }{\mu } \right)^{1/d}, 2^n \epsilon \, \right)  \; ,
  \end{equation}
  
  Now consider $ \sparfom $ with $ \fom = \subgrad $.
  
  From $ K_{\subgrad}(\delta, \bar{\epsilon}  ) :=   (M \delta/ \bar{\epsilon} )^2  $ follows
\[  
    K_{\subgrad} \left( (5 \cdot 2^n \epsilon/\mu)^{1/d}, \, 2^n \epsilon \right)  \leq  \left(  \frac{M (5 \cdot 2^n \epsilon/ \mu )^{1/d}}{2^n \epsilon} \right)^2 = \left( \frac{ 5^{1/d} M }{\mu^{1/d} \epsilon^{1 - 1/d}} \right)^2   \left( \frac{1}{ 4^{1 - 1/d}} \right)^n \; , 
\]
and hence (\ref{eqn.ed})  is bounded above by
\[ 
\bar{N} + 1 + 3 \left( \frac{ 5^{1/d} M }{\mu^{1/d} \epsilon^{1 - 1/d}} \right)^2  \sum_{n=-1}^{\bar{N}} \left( \frac{1}{ 4^{1 - 1/d}} \right)^n \; . 
\] 
The implication for $ d = 1 $  is immediate. On the other hand, since for $ d > 1 $, 
\[  
  \sum_{n=-1}^{\bar{N}} \left( \frac{1}{ 4^{1 - 1/d}} \right)^n < 
 \min \left\{ \sum_{n=-1}^{\infty} \left( \frac{1}{ 4^{1 - 1/d}} \right)^n, \, \, \bar{N} + 5  \right\} 
 =   \min 
  \left\{ \frac{ 16^{1 - 1/d}}{4^{1 - 1/d} - 1}, \, \bar{N} + 5 \right\} \; ,  
\]  
the implication for $ d > 1 $  is immediate, too. 

To obtain upper bounds that hold regardless of whether the inequality $ f( \mathbf{x_0}) - f^* < 5 \cdot 2^N \epsilon $ is satisifed, then according to Theorem~\ref{thm.ea}, simply substitute $ N $ for $ \bar{N} $, and add 
(\ref{eqn.ec}), completing the proof for $ \subgrad $.

In the same manner that above proof for $ \subgrad $  is nearly identical to proof of Corollary~\ref{cor.cc}, the proofs for $ \accel $ and $ \smooth $ are nearly identical to the proofs of Corollaries~\ref{cor.cd}  and \ref{cor.ce}, respectively. The proofs for $ \accel $ and $ \smooth $ are thus left to the reader.  \hfill $ \Box $

\section{{\bf  Proof of Theorem~\ref{thm.ea}}}  \label{sect.f} 

The theorem is obtained through a sequence of results in which we speak of $ \fom_n $ ``updating'' at a point $ x $. The starting point $ \mathbf{x_0} $ is considered to be the first update point for every $ \fom_n $.  After $ \sparfom $ has started, then for $ n < N $, updating at $ x $ means restarting at $ x $. For $ \fom_{N} $, updating at $ x $ is the same as having computed $ x $ satisfying the current task of $ \fom_{N} $, in which case $ x $ becomes the new ``designated point'' and is sent to the inbox of $ \fom_{N-1} $.

\begin{lemma} \label{lem.fa} 
Assume that at time $ t $, $ \fom_n $ updates at $ x $ satisfying $ f(x) - f^* \geq 2 \cdot 2^n \epsilon $. Then no later than time $ t + K_{\fom}( D(f(x)), 2^n \epsilon ) $, $ \fom_n $ updates again.
\end{lemma}
\noindent {\bf Proof:}  Indeed, if $ \fom_n $ has not updated by the specified time, then at that time, it has computed a point $ y $ satisfying $ f(y) - f^* \leq 2^n \epsilon $, simply by definition of the function $ K_{\fom} $, the assumption that $ \fom_n $ performs one iteration in each time period, and the assumption that time periods are of unit length. Since
\[    f(y) = f(x) + (f(y) - f^*) + (f^* - f(x)) \leq f(x) + 2^n \epsilon  - 2 \cdot 2^n \epsilon = f(x) - 2^n \epsilon  \; , \]
in that case an update happens precisely at the specified time. \hfill $ \Box $
 
\begin{prop}  \label{prop.fb}  
Assume that at time $ t $, $ \fom_n $ updates at $ x $ satisfying $ f(x) - f^* \geq 2 \cdot 2^n \epsilon $. Let $ \mathbf{j}  := \lfloor \frac{ f(x) - f^* }{ 2^n \epsilon } \rfloor - 2 $. Then $ \fom_n $ updates at a point $ \bar{x} $ satisfying $ f( \bar{x}) - f^* < 2 \cdot 2^n \epsilon $ no later than time
\begin{equation}   \label{eqn.fa} 
            t + \sum_{j=0}^{\mathbf{j}}   K_{\fom}( D( f(x) - j \cdot 2^n \epsilon), 2^n \epsilon) \; . 
\end{equation} 
\end{prop} 
\noindent {\bf Proof:} Note that $ \mathbf{j} $ is the integer satisfying
\begin{equation}  \label{eqn.fb} 
       f(x) - (\mathbf{j}+1)\cdot 2^n \epsilon < f^* + 2 \cdot 2^n \epsilon \leq    f(x) - \mathbf{j}\cdot 2^n \epsilon \; . 
\end{equation}
        
Lemma \ref{lem.fa}  implies $ \fom_n $ has a sequence of update points $ x_0, x_1, \ldots, x_J, x_{J+1} $, where $ x_0 = x $,  
\begin{gather} f(x_{j+1}) \leq f(x_j) - 2^n \epsilon \quad  \textrm{for all $ j = 0, \ldots J $} \; ,   \label{eqn.fc}  \\ f(x_{J+1}) < f^* + 2 \cdot 2^n \epsilon \leq f(x_J) \; , \label{eqn.fd} 
\end{gather}
and where $ x_{J+1} $ is obtained no later than time
\begin{equation}  \label{eqn.fe} 
      t + \sum_{j=0}^{J} K_{\fom}( D( f(x_j)), 2^n \epsilon) \; . 
 \end{equation} 
 
 By induction, (\ref{eqn.fc}) implies $ f(x_j) \leq f(x) - j \cdot 2^n \epsilon $ ($ j = 0,1, \ldots, J+1 $), which has as a consequence that 
 \[ K_{\fom}(D(f(x_j)),2^n \epsilon) \leq K_{\fom}(D(f(x) - j \cdot 2^n \epsilon),2^n \epsilon) \; ,  \]
and also has as a consequence -- due to the left inequality in (\ref{eqn.fb})  and the right inequality in (\ref{eqn.fd})  -- that $ J \leq \mathbf{j} $.   
Hence, the quantity (\ref{eqn.fe})  is bounded above by (\ref{eqn.fa}), completing the proof for the choice $ \bar{x} = x_{J+1} $ (an appropriate choice for $ \bar{x} $  due to the left inequality in (\ref{eqn.fd})). \hfill $ \Box $
\vspace{2mm}

The following corollary replaces the upper bound (\ref{eqn.fa}) with a bound which depends on quantities $ f^* + i \cdot 2^n \epsilon $ instead of $ f(x) - j \cdot 2^n \epsilon $. This is perhaps the key conceptual step in the proof of Theorem~\ref{thm.ea}.
From this modified bound, Corollary~\ref{cor.fd} shows that not long after $ \fom_n $ satisfies the condition $ f(x) - f^* < 5 \cdot 2^n \epsilon $, $ \fom_{n-1} $ updates at a point $x'$ satisfying the similar condition $ f(x') - f^* < 5 \cdot 2^{n-1} \epsilon $. Inductively applying this result, Corollary~\ref{cor.fe} bounds how long it takes for $ \fom_{-1} $ to find a $(2\cdot 2^{-1}\epsilon)$-optimal solution. This sequence of results is the heart of our analysis and positions us to prove the full runtime bound claimed by Theorem~\ref{thm.ea}.

\begin{cor} \label{cor.fc} 
Assume that at time $ t $, $ \fom_n $ updates at $ x $ satisfying $  f(x) - f^* < \mathbf{i}  \cdot 2^n \epsilon $, 
where $ \mathbf{i}  $ is an integer and $ \mathbf{i} \geq 3 $. Then $ \fom_n $ updates at a point $ \bar{x} $ satisfying  $ f( \bar{x}) - f^* < 2 \cdot 2^n \epsilon $ no later than time
\begin{align}
  & t + \sum_{i=3}^{ \mathbf{i}} K_{\fom}(D_{n,i}, 2^n \epsilon) \label{eqn.ff}  \\
   & \qquad  \qquad  \textrm{where} \, \, D_{n,i} = \min \{ D(f^* + i \cdot 2^n \epsilon), D(f(\mathbf{x_0}  ) ) \} \; . \label{eqn.fg} 
   \end{align}
 \end{cor}
 \noindent {\bf Proof:} Clearly, we may assume $ f(x) \geq f^* + 2 \cdot 2^n \epsilon $. Let $ \mathbf{j} $ be as in Proposition~\ref{prop.fb}, and hence the time bound (\ref{eqn.fa})  applies for updating at some $ \bar{x} $ satisfying $ f( \bar{x}) - f^* < 2 \cdot 2^n \epsilon $.  
  
Note that for all non-negative integers $ j $, 
 \begin{equation}  \label{eqn.fh} 
     f(x) - j \cdot 2^n \epsilon < f^* + (\mathbf{i} - j) \cdot 2^n \epsilon \; . 
     \end{equation}
Consequently, since $ \hat{f} \mapsto K_{\fom}( D( \hat{f}), 2^n \epsilon ) $ is an increasing function, and since $ f(x) \leq f( \mathbf{x_0}) $ (due to $ x $ being an update point),   the quantity (\ref{eqn.fa})  is bounded above by
\begin{equation}  \label{eqn.fi} 
   t + \sum_{i = (\mathbf{i} - \mathbf{j})}^{\mathbf{i}} K_{\fom}(D_{n,i}, 2^n \epsilon) \; . 
\end{equation}
Finally, as $ \mathbf{j} $ is the integer satisfying the inequalities (\ref{eqn.fb}), the rightmost of those inequalities, and (\ref{eqn.fh})  for $ j = \mathbf{j} $, imply $ \mathbf{i} - \mathbf{j} \geq 3 $, and thus (\ref{eqn.fi}) is bounded from above by (\ref{eqn.ff}). \hfill $ \Box $

\begin{cor}  \label{cor.fd} 
Let $  n > -1 $. Assume that at time $ t $, $ \fom_{n} $ updates at $ x $ satisfying $ f(x) - f^* < 5 \cdot 2^n \epsilon $. Then no later than time
\[ 
     t + 1  + \sum_{i=3}^5 K_{\fom}( D_{n,i}, 2^n \epsilon ) 
\] 
(where $ D_{n,i} $ is given by (\ref{eqn.fg})), $ \fom_{n-1} $ updates at $ x' $ satisfying $ f( x')  - f^* < 5 \cdot 2^{n-1} \epsilon $.
\end{cor}
\noindent {\bf Proof:}  By Corollary~\ref{cor.fc}, no later than time 
\[ t +  \sum_{i=3}^5 K_{\fom}( D_{n,i}, 2^n \epsilon )  \; , \]
  $ \fom_n $ updates at $ \bar{x} $ satisfying $ f(\bar{x}) < f^* + 2 \cdot 2^n \epsilon $. When $ \fom_n $ updates at $ \bar{x} $, the point is sent to the inbox of $ \fom_{n-1} $, where it is available to $ \fom_{n-1} $ at the beginning of the next period.

Either $ \fom_{n-1} $ restarts at $ \bar{x} $, or its most recent restart point $ \hat{x}  $ satisfies $ f( \bar{x}  ) > f( \hat{x} ) - 2^{n-1} \epsilon $. In the former case, $ \fom_{n-1} $ restarts at a point $ x' = \bar{x} $    satisfying $  f(x') < f^* + 2 \cdot 2^n \epsilon = f^* + 4 \cdot 2^{n-1} \epsilon $, whereas in the latter case it has already restarted at a point $ x' = \hat{x}  $ satisfying $ f(x') < f( \bar{x}) + 2^{n-1} \epsilon < f^* + 5 \cdot 2^{n-1} \epsilon $. \hfill $ \Box $

\begin{cor}  \label{cor.fe} 
For any $ \bar{N} \in \{ -1, \ldots, N \} $, assume that at time $ t $, $ \fom_{\bar{N}} $ updates at $ x $ satisfying $ f(x) - f^* < 5 \cdot 2^{\bar{N}}  \epsilon $. Then no later than time
\[  
     t + \bar{N}+1 + \sum_{n = -1}^{\bar{N}} \sum_{i=3}^5 K_{\fom}(D_{n,i}, 2^n \epsilon ) \; , 
\] 
$ \sparfom $ has computed an $ \epsilon $-optimal solution.
\end{cor}
\noindent {\bf Proof:}  If $ \bar{N} = -1 $, Corollary \ref{cor.fc} implies the additional time required by $ \fom_{-1} $ to compute a $ (2 \cdot 2^{-1} \epsilon)  $-optimal solution $ \bar{x} $  is bounded from above by
\begin{equation}  \label{eqn.fj} 
      \sum_{i=3}^5 K_{\fom}( D_{-1,i}, \,  2^{-1} \epsilon ) \; . 
      \end{equation} 
The present corollary thus is established for the case $ \bar{N} = -1 $.

On the other hand, if $ \bar{N} > - 1 $, induction using Corollary~\ref{cor.fd}  shows that no later than time
\[  
   t +  \bar{N}+1  + \sum_{n=0}^{\bar{N}} \sum_{i=3}^5 K_{\fom}( D_{n,i}, 2^n \epsilon ) \; , 
\] 
$ \fom_{-1} $ has restarted at $ x $ satisfying $ f(x) - f^* < 5 \cdot 2^{-1} \epsilon $. Then, as above, the additional time required by $ \fom_{-1} $ to compute an $ \epsilon $-optimal solution does not exceed (\ref{eqn.fj}).   \hfill $ \Box $
\vspace{2mm}

We are now in position to prove Theorem \ref{thm.ea}, which we restate as a corollary for the reader's convenience.  Recall that $ \mathbf{x_0} $ is the point at which each $ \fom_n $ is started at time $ t = 0 $. 

\begin{cor}  \label{cor.ff} 
If $ f(\mathbf{x_0}) - f^* < 5 \cdot 2^N \epsilon $, then $ \sparfom $ computes an $ \epsilon $-optimal solution within time
\begin{align} 
  & \bar{N}+ 1 +  3 \sum_{n=-1}^{\bar{N}}  K_{\fom}\left(   D_n , 2^n \epsilon \right)  \label{eqn.fk}   \\   
& \qquad \qquad   \textrm{with } \,    D_n :=   \min \{ D(f^* + 5 \cdot 2^n \epsilon), D( f( \mathbf{x_0})) \}  \; , \nonumber 
\end{align}   
where $ \bar{N} $ is the smallest integer satisfying both $ f(\mathbf{x_0}) - f^* < 5 \cdot 2^{\bar{N}} \epsilon $ and $ \bar{N} \geq -1 $.

In any case, $ \sparfom $ computes an $ \epsilon $-optimal solution within time
\[ 
\mathbf{T}_N + 
  K_{\fom}\big( \dist( \mathbf{x_0}, X^*), 2^{N} \epsilon \big)   \; , 
\]  
 where $ \mathbf{T}_N  $ is the quantity obtained by substituting $ N $ for $ \bar{N} $ in (\ref{eqn.fk}).
\end{cor}
\noindent {\bf Proof:}   For the case that $ f( \mathbf{x_0}) - f^* < 5 \cdot 2^{N} \epsilon $, the time bound (\ref{eqn.fk})  is immediate from Corollary~\ref{cor.fe}  and the fact that $ K_{\fom}(D_{n,i}, 2^n \epsilon ) \leq  K_{\fom}(D_{n,5}, 2^n \epsilon )  $ for $ i = 3,4,5 $.   

In any case, because $ \fom_{N} $ never restarts, within time 
\begin{equation}  \label{eqn.fl} 
        K_{\fom}( \dist( \mathbf{x_0}, X^*), 2^{N}\epsilon) \; , 
        \end{equation}
$ \fom_{N} $ computes a point $ x $ satisfying $ f(x) - f^* \leq 2^{N} \epsilon $. If $ x $ is not an update point for $ \fom_{N} $, then the most recent update point $ \hat{x} $ satisfies $ f(\hat{x}) < f(x) + 2^{N} \epsilon $. Hence, irrespective of whether $ x $ is an update point, by the time $ \fom_{N} $ has computed $ x $, it has obtained an update point $ x' $ satisfying $ f(x') - f^* < 2 \cdot 2^{N} \epsilon < 5 \cdot 2^{N} \epsilon $. Consequently, relying on Corollary~\ref{cor.fe} for $ \bar{N} = N $, as well as on the relations $ K_{\fom}(D_{n,i}, 2^n \epsilon ) \leq  K_{\fom}(D_{n,5}, 2^n \epsilon )  $ for $ i = 3,4,5 $, the time required for $ \sparfom $ to compute an $ \epsilon$-optimal solution does not exceed the value (\ref{eqn.fl})    plus the value (\ref{eqn.fk}), where in (\ref{eqn.fk}),  $ N $ is substituted for $ \bar{N} $. \hfill $ \Box $

\section{{\bf  Asynchronous Parallel Scheme ($ \aparfom $)}} \label{sect.g} 

Motivation for developing an asynchronous scheme arises in two ways. One consideration is that some first order methods require a number of oracle calls that varies significantly among iterations (e.g., adaptive methods utilizing backtracking). If, as in $ \sparfom $, the algorithms $ \fom_n $ are made to wait until each of them has finished an iteration, the time required to obtain an $ \epsilon $-optimal solution might be greatly affected.

The other source of motivation comes from an inspection of the analysis given for $ \sparfom $ in \S\ref{sect.f}. The analysis reveals there potentially is much to be gained from an asynchronous scheme, even for methods which, like $ \subgrad $ and $ \accel $, require the same number of oracle calls at every iteration.  In particular, the gist of Corollary~\ref{cor.fd}  is that the role $ \fom_n $ plays in guiding $ \fom_{n-1} $ is critical only after $ \fom_n $ has first obtained an update point $ x $ satisfying $ f(x) < f^* + 5 \cdot 2^n \epsilon $. Moreover, after $ \fom_n $ has such an update point $ x $, $ \fom_n $ will send $ \fom_{n-1} $ at most four points (and possibly $ x $), simply because $ \fom_n $ sends a point only when the objective has been decreased by at least $ 2^n \epsilon $, and such a decrease can happen at most four times (due to $ f(x) < f^* + 5 \cdot 2^n \epsilon $). Thus, it is only points among the last five that $ \fom_n $ sends to $ \fom_{n-1} $ which play an essential role in the efficiency of $ \sparfom $. 

However, during the time that $ \fom_n $ is performing iterations to arrive at the last few points to be sent to $ \fom_{n-1} $,  an algorithm $ \fom_m $, for $ m \ll n $, can be sending scads of messages to $ \fom_{m-1} $. It is thus advantageous if $ \fom_n $ is disengaged, to the extent possible, from the timeline being followed by $ \fom_m $, especially if all messages go through a single server. Our asynchronous scheme keeps $ \fom_n $ appropriately disengaged. 

\subsection{Details of $ \aparfom $} \label{sect.ga} 

We now turn to precisely specifying the asynchronous scheme. Particular attention has to be given to the delay between when a message is sent by $ \fom_n $ and when it is received by $ \fom_{n-1} $, and what is to be done if $ \fom_{n-1} $ happens to be in the middle of a long iteration when the message arrives. A number of other details also have to be specified in order to prove results regarding the scheme's performance. Specifications different than the ones we describe can lead to similar theoretical results.

We refer to ``epochs'' rather than time periods to denote the duration that $\fom_n $ spends on each of its tasks. Each $ \fom_n $ has its own epochs. 

The first epoch for every $ \fom_n $ begins at time zero, when $ \fom_n $ is started at $ x_{n,0} = \mathbf{x_0} \in Q $, the same point for every $ \fom_n $. 

As for the synchronous scheme, $ \fom_{N} $ never restarts, nor receives messages. Its only epoch is of infinite length.  During the epoch, $ \fom_{N} $ proceeds as before, focused on accomplishing the task of computing an iterate satisfying $ f(x_{N,k}) \leq f(\bar{x}_N ) - 2^{N} \epsilon $ where $ \bar{x}_N $ is the most recent ``designated point.'' If iterate $ x_{N,k} $ satisfies the inequality, it becomes the new designated point, and is sent to the inbox of $ \fom_{N-1} $, the only difference with $ \sparfom $ being that the point arrives in the inbox  no more than $ \ttransit $ units of time after being sent, where $ \ttransit > 0 $.  The same transit time holds for messages sent by any $ \fom_n $ to $ \fom_{n-1} $. (After fully specifying the scheme, we explain that the assumption of a uniform time bound on message delivery is not overly restrictive.)

For $ n < N $, the beginning of a new epoch coincides with $ \fom_n $ obtaining a point $ x $ satisfying the inequality of its task, that is, $ f(x) \leq f(x_{n,0}) - 2^n \epsilon $ where $ x_{n,0} $ is the most recent (re)start point for $ \fom_n $. 

As for the synchronous scheme, for $ n < N $, $ \fom_n $ can receive messages from $ \fom_{n+1} $. Now, however, it is assumed that when a message arrives in the inbox, 
a ``pause'' for $ \fom_n $ begins, lasting a positive amount of time, but no more than $ \tpause $ time units. The pause is meant to capture various delays that can happen upon receipt of a new message, including time for $ \fom_n $ to come to a good stopping place in its calculations before actually going to the inbox.

If another message from $ \fom_{n+1} $ arrives in the inbox during a pause, then immediately the old message is overwritten by the new one, the old pause is cancelled and a new pause begins. 

If during an epoch for $ \fom_n $, a message arrives after $ \fom_n $ has restarted (i.e., after $ \fom_n $ has begun making iterations), then during the pause, $ \fom_n $ determines whether the point $ x' $ sent by $ \fom_{n+1} $ satisfies the inequality of $ \fom_n $'s current task.  If during the pause, another message from $ \fom_{n+1} $ arrives, overwriting $ x' $ with a new point $ x'' $, then in the new pause, $ \fom_n $ turns attention to determining if $ x'' $ satisfies the inequality (indeed, $ x'' $  is a better candidate than $ x' $  because $ f(x'') \leq  f(x') - 2^{n+1} \epsilon $). And so on, until the sequence of contiguous pauses ends.\footnote{The sequence of contiguous pauses is finite, because (1) $ \fom_{n+1} $ sends a message only when it has obtained $ x $ satisfying $ f(x) \leq f(x_{n+1,0}) - 2^{n+1} \epsilon $, and (2) we assume $ f^* $ is finite.} Let $ x $ be the final point.

If $ x $ satisfies the inequality of the task for $ \fom_n $, then immediately at the end of the (contiguous) pause, a new epoch begins for $ \fom_n $. On the other hand, if $ x $ fails to satisfy the inequality, then the current epoch continues, with $ \fom_n $ returning to making iterations.

The other manner an epoch for $ \fom_n $ can end is by $ \fom_n $ itself computing an iterate that satisfies the inequality. At the instant that $ \fom_n $ computes such an iterate, a new epoch begins.

It remains to describe what occurs at the beginning of an epoch for $ \fom_n $. 

As already mentioned, at time zero -- the beginning of the first epoch for all $ \fom_n $ -- every $ \fom_n $ is started at $ x_{n,0} = \mathbf{x_0} $.  No messages are sent at time zero.   

For $ n < N $, as stated above, the beginning of a new epoch occurs precisely when $ \fom_n $ has obtained a point $ x $ satisfying the inequality of its task.

If no message from $ \fom_{n+1} $ arrives simultaneously with the beginning of a new epoch for $ \fom_n $, then $ \fom_n $ instantly\footnote{We assume the listed chores are accomplished instantly by $ \fom_n $. Assuming positive time is required would have negligible effect on the complexity results, but would make notation more cumbersome.}    (1) relabels $ x $ as $ x_{n,0} $ and updates its task accordingly, (2) restarts at $ x_{n,0} $, and (3) sends a copy of $ x_{n,0}  $ to the inbox of $ \fom_{n-1} $ (assuming $ n > - 1 $), where it will arrive no later than time $ t + \ttransit $, with $ t $ being the time the new epoch for $ \fom_n $ has begun.

On the other hand, if a message from $ \fom_{n+1} $ arrives in the inbox exactly at the beginning of the new epoch, then $ \fom_n $ immediately pauses. Here, $ \fom_n $ makes different use of the pause than what occurs for pauses happening after $ \fom_n $ has restarted. Specifically, letting $ x' $ be the point in the inbox, then during the pause, $ \fom_n $ determines whether $ f(x') < f(x) $ -- if so, $ x' $ is preferred to $ x $ as the restart point. If during the pause, another message arrives, overwriting $ x' $ with $ x'' $, then during the new pause, $ \fom_n $ determines whether $ x'' $ is preferred to $ x  $ ($ x'' $ is definitely preferred to $ x' $ -- indeed, $ f(x'') \leq  f(x') - 2^{n+1} \epsilon $). And so on, until the sequence of contiguous pauses ends, at which time $ \fom_n $ instantly (1) labels the most preferred point as $ x_{n,0} $ and updates its task accordingly, (2) restarts at $ x_{n,0} $, and (3) sends a copy of $ x_{n,0} $ to the inbox of $ \fom_{n-1} $ (assuming $ n > - 1 $). 

The description of $ \aparfom $ is now complete.

\subsection{Remarks} \label{sect.gb} 
\begin{itemize}

\item  Regarding the transit time, a larger value for $ \ttransit $ might be chosen to reflect the possibility of greater congestion in message passing, as could occur if all messages go through a queue on a single server. We emphasize, however, that a long queue of unsent messages can be avoided, because (1) if two messages from $ \fom_n $ are queued, the later message contains a point that has better objective value than the earlier message, and (2) as indicated above, only the last few messages (at most four, perhaps none) sent by $ \fom_n $ play an essential role.  Thus, to avoid congestion, when a second message from $ \fom_n $ arrives in the queue, delete the earlier message. So long as this policy of deletion is also applied to all of the other copies $ \fom_m $, the fact that only the last few messages from $ \fom_n $ are significant then implies there is not even a need to move the second message from $ \fom_n $ forward into the place that had been occupied by the deleted message -- any truly essential message from $ \fom_n $ will naturally move from the back to the front of the queue within time proportional to $ N+1 $ after its arrival.  Managing the queue in this manner -- deleting a message from $ \fom_n $ as soon as another message from $ \fom_n $ arrives in the queue -- ensures that $ \ttransit $ can be chosen proportional to $ N+1 $ in the worst case of there being only a single server devoted to message passing.
\item  As for the synchronous scheme, the ideal arrangement is, of course, for each $ n > -1 $, to have apparatus dedicated solely to the sending of messages from $ \fom_n $ to $ \fom_{n-1} $. Then $ \ttransit $ can be chosen as a constant independent of the size of $ N $.
\item For many first-order methods, it is natural to choose $ \tpause $ to be the maximum possible number of oracle calls in an iteration, in which case a pause for $ \fom_n $ can be interpreted as the time needed for $ \fom_n $ to complete its current iteration before checking whether a message is in its inbox. 
\end{itemize}

The synchronous scheme can be viewed as a special case of the asynchronous scheme by choosing $ \ttransit = 1 $ and $ \tpause = 0 $ (rather, can be viewed as a limiting case, because for $ \aparfom $, the length of a pause is assumed to be positive).   Indeed, choices in designing the asynchronous scheme were made with this in mind, the intent being that the main result for the synchronous scheme would be a corollary of a theorem for the asynchronous scheme.  However, as the foremost ideas for both schemes are the same, and since a proof for the synchronous scheme happens to have greater transparency, in the end we chose to first present the theorem and proof for the synchronous scheme.

\section{{\bf  Theory for $ \aparfom $}}  \label{sect.h}

Here we state the main theorem regarding $ \aparfom $, and deduce from it a corollary for Nesterov's universal fast gradient method \cite{nesterov2015universal}. 

\subsection{Main theorem for $ \aparfom $} \label{sect.ha} 
So as to account for the possibility of variability in the amount of work among iterations, we now rely on a function $ T_{\fom}(\delta, \bar{\epsilon})  $ assumed to provide an upper bound on the amount of time (often proportional to the total number of oracle calls) required by $ \fom $ to compute an $ \bar{\epsilon}  $-optimal solution when $ \fom $ is initiated at an arbitrary point $ x_0 \in Q $ satisfying $ \dist(x_0, X^*) \leq \delta $.

Recall that for scalars $ \hat{f} \geq f^* $, 
\[  
   D(\hat{f}) :=  \sup \{ \dist(x,X^*) \mid  x \in Q  \textrm{ and } f(x) \leq \hat{f} \} \; . 
\]

\begin{thm} \label{thm.ha}  
If $ f(\mathbf{x_0}) - f^* < 5 \cdot 2^N \epsilon $, then $ \aparfom $ computes an $ \epsilon $-optimal solution within time
\begin{align} 
  & (\bar{N}+1) \ttransit + 2(\bar{N}+2) \tpause +  3 \sum_{n=-1}^{\bar{N}}  T_{\fom}\left(   D_n , 2^n \epsilon \right)  \label{eqn.ha}   \\   
& \qquad \qquad   \textrm{with } \,    D_n :=   \min \{ D(f^* + 5 \cdot 2^n \epsilon), D( f( \mathbf{x_0})) \}  \; , \nonumber 
\end{align}   
where $ \bar{N} $ is the smallest integer satisfying both $ f(\mathbf{x_0}) - f^* < 5 \cdot 2^{\bar{N}} \epsilon $ and $ \bar{N} \geq -1 $.

In any case, $ \aparfom $ computes an $ \epsilon $-optimal solution within time
\[ 
\mathbf{T}_N + 
  T_{\fom}\big( \dist( \mathbf{x_0}, X^*), 2^{N} \epsilon \big)   \; , 
\]  
 where $ \mathbf{T}_N  $ is the quantity obtained by substituting $ N $ for $ \bar{N} $ in (\ref{eqn.ha}).
\end{thm}

The proof is deferred to Appendix \ref{sect.l}, due to its similarities with the proof of Theorem~\ref{thm.ea}.
\vspace{1mm}

\noindent {\bf Remarks:} 
\begin{itemize}

\item As observed in \S\ref{sect.g}, a larger value for $ \ttransit $ might be chosen to reflect the possibility of greater congestion in message passing, although by appropriately managing communication, in the worst case (in which all communication goes through a single server), a message from $ \fom_n $ would reach $ \fom_{n-1} $ within time proportional to $ N+1 $.  For $ \ttransit $ being proportional to $ N+1 $, we see from (\ref{eqn.ha})  that the impact of congestion on $ \aparfom $ is modest, adding an amount of time proportional to $ (N+1)^2 $.  Compare this with the synchronous scheme, where to incorporate congestion in message passing, every time period would be lengthened, in the worst case to length $ 1 + \kappa \cdot   (N+1) $ for some postive constant $ \kappa  $. The time bound (\ref{eqn.ea})  would then be multiplied by $ 1 + \kappa \cdot   (N+1) $, quite different than adding only a single term proportional to $ (N+1)^2 $.
\item Unlike the synchronous scheme, for the asynchronous scheme we know of no sequential analogue that in general is truly natural. 
\end{itemize}

\subsection{A corollary} \label{sect.hb} 
To provide a representative application of Theorem \ref{thm.ha}, we consider Nesterov's universal fast gradient method \cite[\S4]{nesterov2015universal}-- denoted $ \univ $ -- which applies whenever $ f $ has H\"{o}lder continous gradient with exponent $ \nu $ ($ 0 \leq \nu \leq 1 $), meaning
\[  M_{\nu} := \sup \left\{ \smfrac{ \| \grad f(x) - \grad f(y) \| }{ \| x - y \|^{ \nu}} \mid x,y \in Q, \, x \neq y \right\} \,  < \, \infty \quad \textrm{(i.e., is finite)} \; .  \]
(If $ \nu = 0 $ then $ f $ is $ M_0 $-Lipschitz on $ Q $, whereas if $ \nu = 1 $, $ f $ is $ M_1 $-smooth on $ Q $.)

The input to $ \univ $  consists of the desired accuracy $ \bar{\epsilon}  $, an initial point $ x_0 \in Q $, and a value $ L_0  > 0 $ meant, roughly speaking, as a guess of $ M_{\nu} $. 
Nesterov \cite[(4.5)]{nesterov2015universal} showed the function 
\begin{equation}  \label{eqn.hb} 
   K_{\univ}(\delta, \bar{\epsilon} ) =  4 \left(  M_{\nu} \delta^{1 + \nu}/ \bar{\epsilon}  \right)^{\frac{2}{1 + 3 \nu}} 
   \end{equation} 
   provides an upper bound on the number of iterations sufficient to compute an $ \bar{\epsilon}$-optimal solution if $ \dist(x_0, X^*) \leq \delta $ (where in \cite[(4.5)]{nesterov2015universal} we have substituted $ \xi(x_0,x^*) = \frac{1}{2} \| x_0 - x^* \|^2 $ and used $ 2^{ \frac{2(1 + \nu)}{1 + 3 \nu }} \leq 4 $ when $ 0 \leq \nu \leq 1 $).
    
The number of oracle calls in some iterations, however, can significantly exceed the number in other iterations. Indeed, the proofs in \cite{nesterov2015universal} leave open the possibility that the number of oracle calls made only in the $ k^{th} $ iteration might exceed $ k $. Thus, the scheme $ \aparfom $ is highly preferred  to $ \sparfom $ when $ \fom $ is chosen to be $ \univ $.

For the universal fast gradient method, the upper bound established in \cite{nesterov2015universal} on the number of oracle calls in the first $ k  $ iterations is 
\begin{equation}   \label{eqn.hc} 
   4( k +1) + \log_2 \left( \delta^{\frac{2(1 - \nu)}{1 + 3 \nu}} (1/ \bar{\epsilon} )^{\frac{3(1-\nu)}{1 + 3 \nu}} M_{\nu}^{ \frac{4}{1 + 3 \nu}} \right) - 2 \log_2 L_0 \; , 
\end{equation} 
assuming $ \dist(x_0,X^*) \leq \delta $.
This upper bound depends on an assumption such as 
\begin{equation}  \label{eqn.hd} 
  L_0 \leq \left( \smfrac{1 - \nu }{1+ \nu } \cdot \smfrac{1}{\bar{\epsilon }} \right)^{\smfrac{1 - \nu }{1+ \nu }} M_{ \nu}^{\smfrac{2}{1+ \nu }}  \;  
  \end{equation} 
($ L_0 \leq M_1 $ when $ \nu = 1 $).
Presumably the assumption can be removed by a slight extension of the analysis, resulting in a time bound that differs insignificantly, just as the assumption can readily be removed for the first universal method introduced by Nesterov in \cite{nesterov2013gradient}  (essentially the same algorithm as $ \univ $ when $ \nu = 1 $).\footnote{For that setting, see \cite[Appendix A]{renegar2017accelerated}  for a slight extension to Nesterov's arguments that suffice to remove the assumption.}  
However, as the focus of the present paper is on the simplicity of $ \parfom $ and not on $ \univ $ per se, we do not digress to attempt removing the assumption, but instead make assumptions that ensure the results from \cite{nesterov2015universal}  apply. 

In particular, we assume $ L_0 $ is a positive constant satisfying (\ref{eqn.hd})  for $ \bar{\epsilon} = 2^{-1} \epsilon $, and thus satisfies (\ref{eqn.hd})  when $ \bar{\epsilon} = 2^n \epsilon $ for any $ n \geq -1 $.  Moreover, we assume that whenever $  \univ_n $ ($ n = -1, 0, \ldots N $) is (re)started, the input consists of $ \bar{\epsilon} = 2^n \epsilon $, $ x_{n,0} $ (the (re)start point), and $ L_0 $ (the same value at every (re)start). 

Relying on the same value $ L_0 $ at every (re)start likely results in complexity bounds that are slightly worse in some cases (specifically, when $ d = 1 + \nu $ and $ 0 \leq \nu < 1 $), but still not far from being optimal, as we will see.

In view of (\ref{eqn.hc}), and assuming (\ref{eqn.hd}) holds, it is natural to define
\begin{equation}  \label{eqn.he}
  T_{\univ}(\delta, \bar{\epsilon} ) = 4( K_{\univ}( \delta, \bar{\epsilon} ) + 1) + \log_2 \left( \delta^{\frac{2(1 - \nu)}{1 + 3 \nu}} (1/ \bar{\epsilon} )^{\frac{3(1-\nu)}{1 + 3 \nu}} M_{\nu}^{ \frac{4}{1 + 3 \nu}} \right) - 2 \log_2 L_0\; ,
\end{equation}  
an upper bound on the time (total number of oracle calls) sufficient for $ \univ $ to compute an $ \bar{\epsilon} $-optimal solution when initiated at an arbitrary point $ x_0 \in Q $ satisfying $ \dist(x_0, X^*)\leq\delta $.    In order to reduce notation, for $ \epsilon > 0 $ let
\[  
     {\mathcal C} (\delta, \epsilon) :=  4 + \log_2 \left( \delta^{\frac{2(1 - \nu)}{1 + 3 \nu}} (2/ \epsilon   )^{\frac{3(1-\nu)}{1 + 3 \nu}} M_{\nu}^{ \frac{4}{1 + 3 \nu}} \right) - 2 \log_2 L_0 \; , 
\] 
in which case for $ n \geq -1 $, from (\ref{eqn.he})  and (\ref{eqn.hb})  follows
\begin{align}  
 T_{\univ}(\delta, 2^n \epsilon ) & \leq 4   K_{\univ}( \delta, 2^n \epsilon ) + {\mathcal C} (\delta, \epsilon) \nonumber \\
&  \leq 16  \left( M_{\nu} \delta^{1 + \nu}/ (2^n \epsilon)  \right)^{\frac{2}{1 + 3 \nu}} + {\mathcal C} (\delta, \epsilon) \; . \label{eqn.hf} 
\end{align}
Note that $ {\mathcal C} (\delta, \epsilon) $ is a constant independent of $ \epsilon $ if $ \nu = 1 $, and otherwise grows like $ \log(1/\epsilon) $ as $ \epsilon \rightarrow 0 $.

We assume $ f $ has H\"{o}lderian growth, that is, there exist constants $ \mu > 0 $ and $ d \geq 1 $ for which
\[ 
   x \in Q \textrm{ and } f(x) \leq f( \mathbf{x_0})  \quad \Rightarrow \quad f(x) - f^* \geq \mu \,  \dist(x,X^*)^d \; . 
\] 
Consequently, the values $ D_n $ in Theorem \ref{thm.ha}  satisfy
\begin{equation}  \label{eqn.hg} 
               D_n \leq \min \{ (5 \cdot 2^n \epsilon/\mu)^{1/d}, D(f( \mathbf{x_0})) \} \; .  
\end{equation}
  
For a function $ f $ which has H\"{o}lderian growth and has H\"{o}lder continuous gradient, necessarily the values $ d $ and $ \nu $ satisfy $ d \geq 1 + \nu $ (see \cite[\S1.3]{roulet2017sharpness}).

\begin{cor}  \label{cor.hb} 
Consider $ \aparfom $ with $ \fom = \univ $.  \\
If $ f(\mathbf{x_0}) - f^* < 5 \cdot 2^{N} \epsilon $, then $ \aparfom $ computes an $ \epsilon$-optimal solution in time $ T $ for which 
\begin{align}
  d & = 1 + \nu   \, \,   \Rightarrow \, \,   T \leq (\bar{N}+1) \ttransit  + 2(\bar{N}+2) \tpause + 3(\bar{N}+2) \, {\mathcal C} ( D(f(\mathbf{x_0})), \epsilon) \nonumber \\
&\qquad  \qquad  \qquad  \qquad  \qquad  \qquad  \qquad  \qquad  \qquad  \qquad  \quad  \, \, 
      +  48  (\bar{N} + 2)  \left( 5 M_{\nu}/ \mu \right)^{\frac{2}{1 + 3 \nu}}  \; ,  \label{eqn.hh}  \\
  d &>  1 + \nu \, \,  \Rightarrow \, \,  T \leq (\bar{N}+1) \ttransit  + 2(\bar{N}+2) \tpause + 3(\bar{N}+2) \, {\mathcal C} ( D(f(\mathbf{x_0})), \epsilon)  \nonumber \\
& \quad \qquad  \qquad  \qquad \qquad  +  48 \left( \frac{M_{\nu} (5  / \mu)^{\frac{1 + \nu}{d}}}{\epsilon^{1 - \frac{1 + \nu}{d}}}\right)^{\frac{2}{1 + 3 \nu}} 
 \min \left\{ \frac{4^{(1 - \frac{1 + \nu}{d}) \frac{2}{1 + 3 \nu}}}{2^{(1 - \frac{1 + \nu}{d})\frac{2}{1 + 3 \nu} } - 1}, \bar{N} + 5 \right\}, \label{eqn.hi} 
  \end{align}
where  $ \bar{N} $ is the smallest integer satisfying both $ f( \mathbf{x_0}) < f^* + 5 \cdot 2^{\bar{N}} \epsilon $ and $ \bar{N} \geq -1 $. 

In any case, a time bound is obtained by substituting $ N $ for $ \bar{N} $ above, and adding 
\begin{equation}  \label{eqn.hj} 
16 \left( M_{\nu} \dist( \mathbf{x_0}, X^*)^{1 + \nu}/ (2^N \epsilon) \right)^{\frac{2}{1 + 3 \nu}} + {\mathcal C} ( \dist(\mathbf{x_0},X^*), 2^N \epsilon) \; . 
 \end{equation}    
    \end{cor} 
    \vspace{2mm}

\noindent {\bf Remarks:}  
According to Nemirovski and Nesterov \cite[page 6]{NemNes85} (who state lower bounds when $ 0 < \nu \leq 1 $), for   the easily computable choice $ N = \max \{ 0, \lceil \log_2(1/ \epsilon) \rceil \} $, the time bounds of the corollary would be optimal  -- with regards to $ \epsilon $, $ \mu $, $ d $ and $ M_{ \nu} $  -- for a sequential algorithm, except in the cases when both $ d = 1 + \nu $ and $ 0 < \nu < 1 $, where due to the term ``$ 3( \bar{N}+2) \, {\mathcal C} (D(f( \mathbf{x_0})), \epsilon )  $,'' the bound (\ref{eqn.fh})  grows like $ \log(1/\epsilon)^2 $ as $ \epsilon \rightarrow 0 $, rather than growing like $ \log(1/\epsilon) $. Consequently, in those cases, the total amount of work is within a multiple of $ \log(1/\epsilon)^2 $ of being optimal, whereas in the other cases for which they state lower bounds, the total amount of work is within a multiple of $ \log(1/\epsilon) $ of being optimal.
\vspace{2mm}

 \noindent {\bf Proof of Corollary \ref{cor.hb}:} Define $ D_n $ as in Theorem \ref{thm.ha}, that is, $ D_n = \min \{ D(f^* + 5 \cdot 2^n \epsilon), D( f( \mathbf{x_0})) \} $.  Let $ {\mathcal C} =  {\mathcal C}  ( D(f(\mathbf{x_0}), \epsilon) $. 
 
 Assume $ f( \mathbf{x_0}) - f^* < 5 \cdot 2^N \epsilon $.   By (\ref{eqn.hf})  and (\ref{eqn.hg}),  for $ n \geq -1 $ we have
\begin{align}
  T_{\univ}(D_n, 2^n \epsilon) & \leq {\mathcal C}  +  16 \left( \frac{M_{\nu} (5 \cdot 2^n \epsilon / \mu)^{\frac{1 + \nu}{d}}}{2^n \epsilon} \right)^{\frac{2}{1 + 3 \nu}} \nonumber \\
  & = {\mathcal C}  +  16 \left( \frac{M_{\nu} (5  / \mu)^{\frac{1 + \nu}{d}}}{\epsilon^{1 - \frac{1 + \nu}{d}}}\right)^{\frac{2}{1 + 3 \nu}} \left( \frac{1}{2^{(1 - \frac{1 + \nu}{d})\frac{2}{1 + 3 \nu} }} \right)^n  \; .  \label{eqn.hk} 
 \end{align}
Substituting this into the bound (\ref{eqn.ha})  of Theorem \ref{thm.ha}  establishes the implication (\ref{eqn.hh}) 

For $ d > 1 + \nu $, observe
\begin{align}    
   \sum_{n=-1}^{\bar{N}} \left( \frac{1}{2^{(1 - \frac{1 + \nu}{d})\frac{2}{1 + 3 \nu} }} \right)^n & < \min \left\{ \sum_{n=-1}^{\infty } \left( \frac{1}{2^{(1 - \frac{1 + \nu}{d})\frac{2}{1 + 3 \nu} }} \right)^n, \, \bar{N} + 5 \right\} \nonumber \\
 & = 
   \min \left\{ \frac{4^{(1 - \frac{1 + \nu}{d}) \frac{2}{1 + 3 \nu}}}{2^{(1 - \frac{1 + \nu}{d})\frac{2}{1 + 3 \nu} } - 1}, \, \bar{N} + 5 \right\}
            \; , \label{eqn.hl} 
\end{align}
where for the inequality we have used $ 1 \leq 2^{(1 - \frac{1 + \nu}{d})\frac{2}{1 + 3 \nu} } \leq 4 $.
Substituting (\ref{eqn.hk}) into the bound (\ref{eqn.ha}) of Theorem~\ref{thm.ha}, and then substituting (\ref{eqn.hl}) for the resulting summation, establishes the implication (\ref{eqn.hi}), concluding the proof in the case that $ f( \mathbf{x_0}) - f^* < 5 \cdot 2^n \epsilon $. 

To obtain time bounds when $   f( \mathbf{x_0}) - f^* \geq  5 \cdot 2^{N} \epsilon $, then according to Theorem~\ref{thm.ha}, simply substitute $ N $ for $ \bar{N} $ in the bounds above, and add 
\[  
  T_{\univ}\big( \dist( \mathbf{x_0}, X^*), 2^N \epsilon \big) \; , \]
which due to (\ref{eqn.hf}), is bounded above by (\ref{eqn.hj}). \hfill $ \Box $
\vspace{2mm}

\section{{\bf Numerical Experiments}}  \label{sect.i}  
In this section, we present numerical experiments giving insight into the behavior of the proposed parallel restarting schemes ($\sparfom$ and $\aparfom$) utilizing either $\mathtt{subgrad}$, $\mathtt{accel}$, or $\mathtt{smooth}$. Implementations of both of these methods are available at \footnote{\url{https://github.com/bgrimmer/Parallel-Restarting-Scheme}} in a Jupyter notebook (implemented in Julia). 

\subsection{Restarting $\mathtt{subgrad}$ and $\mathtt{smooth}$ for piecewise linear optimization.}
Our first experiments consider minimizing a generic piecewise linear convex function
\begin{equation}\label{eq:pl-obj}
\min_{x\in\mathbb{R}^n} \max \{ a_i^T x - b_i \mid i = 1, \ldots, m \}  \quad \textrm{where  }  a_i \in \mathbb{R}^n,  b_i \in \mathbb{R} \; .
\end{equation}
We use a synthetic problem instance with $m=2000$ and $n=100$ constructed by sampling each $a_i\in\mathbb{R}^{n}$ from a standard Gaussian distribution and each $b_i\in\mathbb{R}$ from a unit Poisson distribution. In the following subsections, we examine the performance of restarting $\mathtt{subgrad}$, $\mathtt{smooth}$, and $\mathtt{accel}$.
We set a target accuracy of $\epsilon=0.002$, $x_0=(1,\dots,1)$, and $N=14$, which results in using $16$ instances of the given first-order method.

First we apply the subgradient method $\mathtt{subgrad}$ to this piecewise linear problem. Notice that the definition of $\mathtt{subgrad}$ in~\eqref{eqn.ab} is only parameterized by $\epsilon$. Hence applying $\sparfom$ here does not require any form of parameter guessing or tuning to run.

\begin{figure}
	\centering
	\includegraphics[width=\textwidth]{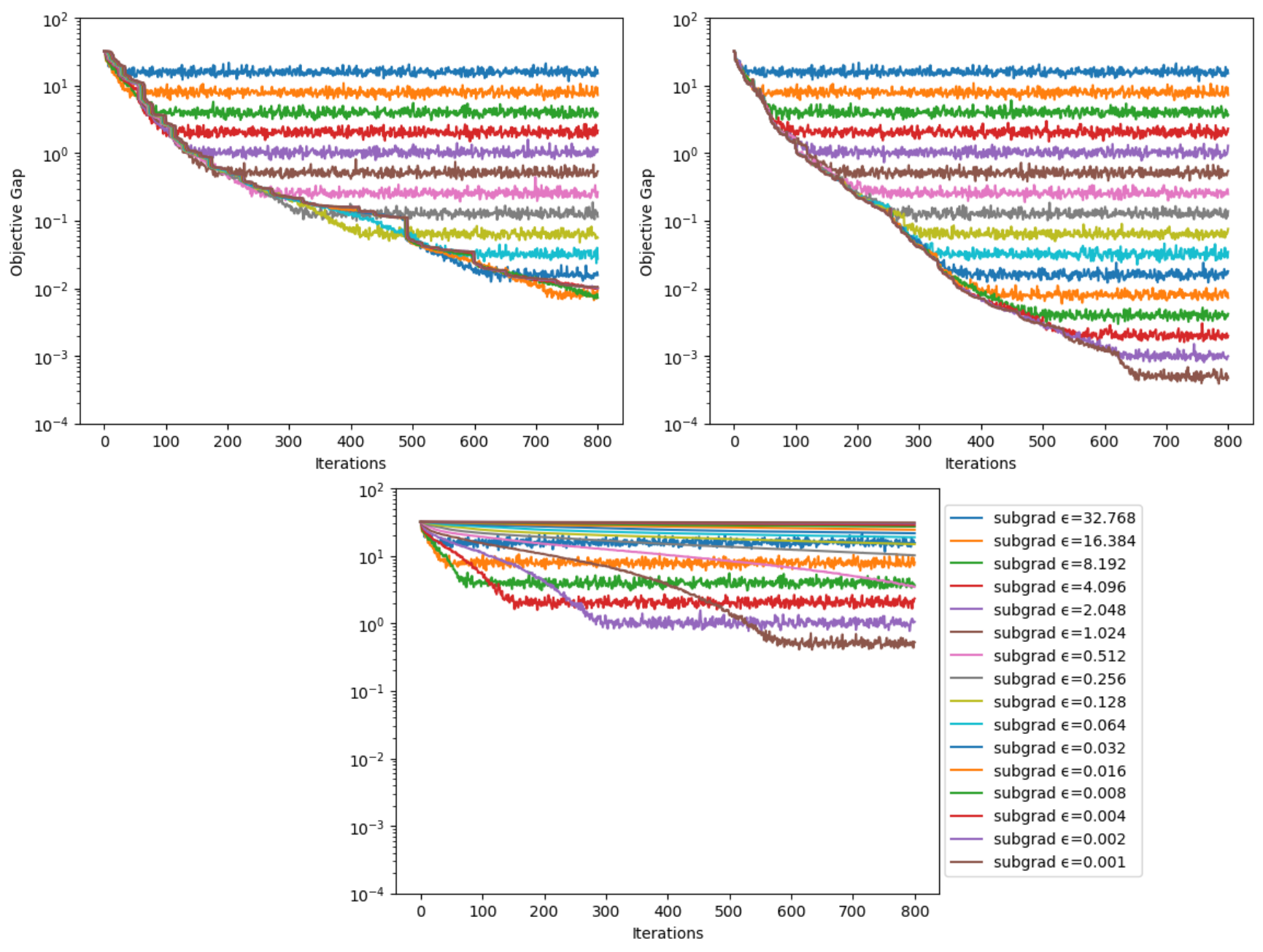}
	
	\caption{The top left plot shows the minimum objective gap of~\eqref{eq:pl-obj} seen by each $\mathtt{subgrad}$ instance in $\sparfom$ with target accuracies $\epsilon=0.001\times 2^k$. The top right plot shows improved convergence of $\sparfom$ when $ \fom_n $ broadcasts its restarts to everyone, not just to $ \fom_{n-1} $. The bottom plot shows the slower convergence of these methods without restarting.} \label{fig:subgrad}
\end{figure}

We apply two variations of $\sparfom$: first as defined in Section~\ref{sect.d} where each instance sends messages with its current iterate to consider for restarting to the next instance whenever it restarts, and second where each instance sends messages to all other instances each iteration. Figure~\ref{fig:subgrad} shows the performance of both variation of $\sparfom$ (in the top left and top right plots respectively) over the course of $800$ iterations as well as the performance (in the bottom plot) of $\mathtt{subgrad}$ without restarting applied with the same values of $\epsilon$ used throughout $\sparfom$.

We see that each $\mathtt{subgrad}$ instance (whether run with $\sparfom$ or alone) flattens out at an objective gap less than its target objective gap of $0.001\times 2^k$. 
Each instance of the subgradient method converges to its target accuracy faster when applied in the restarting scheme, most notably when $\epsilon$ is small. Hence applying our restarting scheme to the subgradient method gives better results than any tuning of the stepsize (controlled by choosing $\epsilon$) could.
Moreover, the $\sparfom$ variation having every instance message every instance further improves the restarting schemes accuracy reached by another order of magnitude.

This piecewise linear problem can be solved more efficiently via the smoothing discussed in Section~\ref{sect.ba}. For any $\eta>0$, consider
\begin{equation}\label{eq:pl-smoothing}
\min_{x\in\mathbb{R}^n} \eta \ln \left( \sum_{i=1}^m \exp ( (a_i^T x - b_i)/ \eta) \right)  \; .
\end{equation}
This a $(\alpha,\beta)=(\max_{i,j}\{a_{ij}^2\}, \ln (m)) $-smoothing of~\eqref{eq:pl-obj} (see~\cite{beck2012smoothing}).
Then $\mathtt{smooth}(\epsilon)$ solves~\eqref{eq:pl-obj} by applying the accelerated method to~\eqref{eq:pl-smoothing} with $\eta=\epsilon/3\beta$, which has a $\alpha/\eta$-smooth objective. Since $\alpha$ and $\beta$ are both known, no parameter guessing or tuning is required to run $\sparfom$ with $\mathtt{smooth}$.

\begin{figure}
	\centering
	\includegraphics[width=\textwidth]{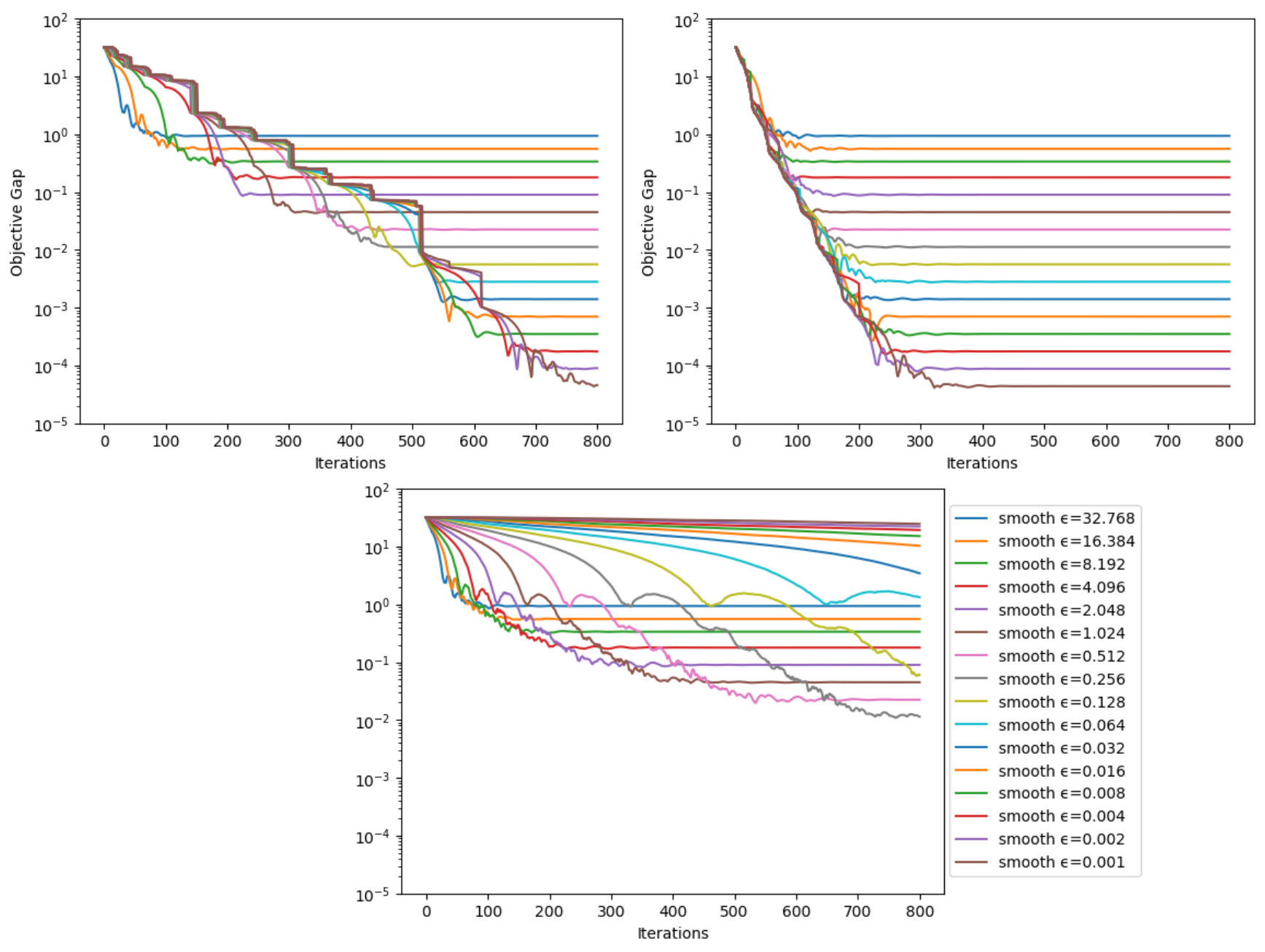}
	
	\caption{The top left plot shows the minimum objective gap of~\eqref{eq:pl-obj} seen by each $\mathtt{smooth}$ instance in $\sparfom$ with target accuracies $\epsilon=0.001\times 2^k$. The top right plot shows improved convergence of $\sparfom$ when $ \fom_n $ broadcasts its restarts to everyone, not just to $ \fom_{n-1} $.  The bottom plot shows the slower convergence of these methods without restarting.} \label{fig:smooth}
\end{figure}

Figure~\ref{fig:smooth} shows the objective gap of each $\mathtt{smooth}$ instance employed by $\sparfom$ (as defined in Section~\ref{sect.d} in the top left plot and with additional message passing in the top right) and the objective gap of each $\mathtt{smooth}$ instance without restarting (in the bottom plot). Much like the subgradient method, each $\mathtt{smooth}$ instance in our scheme strictly dominates the performance of its counterpart without restarting and adding additional message passing to $\sparfom$ notably speeds up the convergence. (That $ .0001 $ accuracy is achieved -- whereas the target is $ .002 $ -- is due to the choice $ \eta = \epsilon/3 \beta $, which always guarantees an $ \epsilon $-optimal solution will be reached, but for this numerical example happened to lead to even greater accuracy.)

\subsection{Restarting $\mathtt{accel}$ for least squares optimization.}
Now consider solving a standard least squares regression problem defined by
\begin{equation}\label{eq:regress}
\min_{x\in\mathbb{R}^n}\ \frac{1}{2m}\|Ax-b\|^2_2
\end{equation}
for some $A\in\mathbb{R}^{m\times n}$ and $b\in\mathbb{R}^m$.
We sample the $A$ and $x^*\in\mathbb{R}^n$ from standard Gaussian distributions with $m=2000$ and $n=1000$ and then set $b = Ax^*$. We initialize our experiments with $x_0=(0,\dots,0)$, $N=30$ and a target accuracy of $\epsilon=10^{-9}$.

We solve this problem with both $\sparfom$ and $\aparfom$ using $32$ parallel versions of $\mathtt{accel}$ (and without any additional message passing). This problem is known to be smooth with a Lipschitz constant $L$ given by the maximum eigenvalue of $A^TA$. Moreover, this problem satisfies quadratic growth (H\"older growth with $d=2$) with coefficient $\mu$ given by the minimum eigenvalue of $A^TA$. Hence this setting is amenable to applying our restarting scheme and should have the convergence rate improved from the sublinear convergence rate of $O(\sqrt{L/\epsilon})$ to the linear convergence rate of $O(\sqrt{L/\mu}\log(1/\epsilon))$.

\begin{figure}
	\centering
	\includegraphics[width=\textwidth]{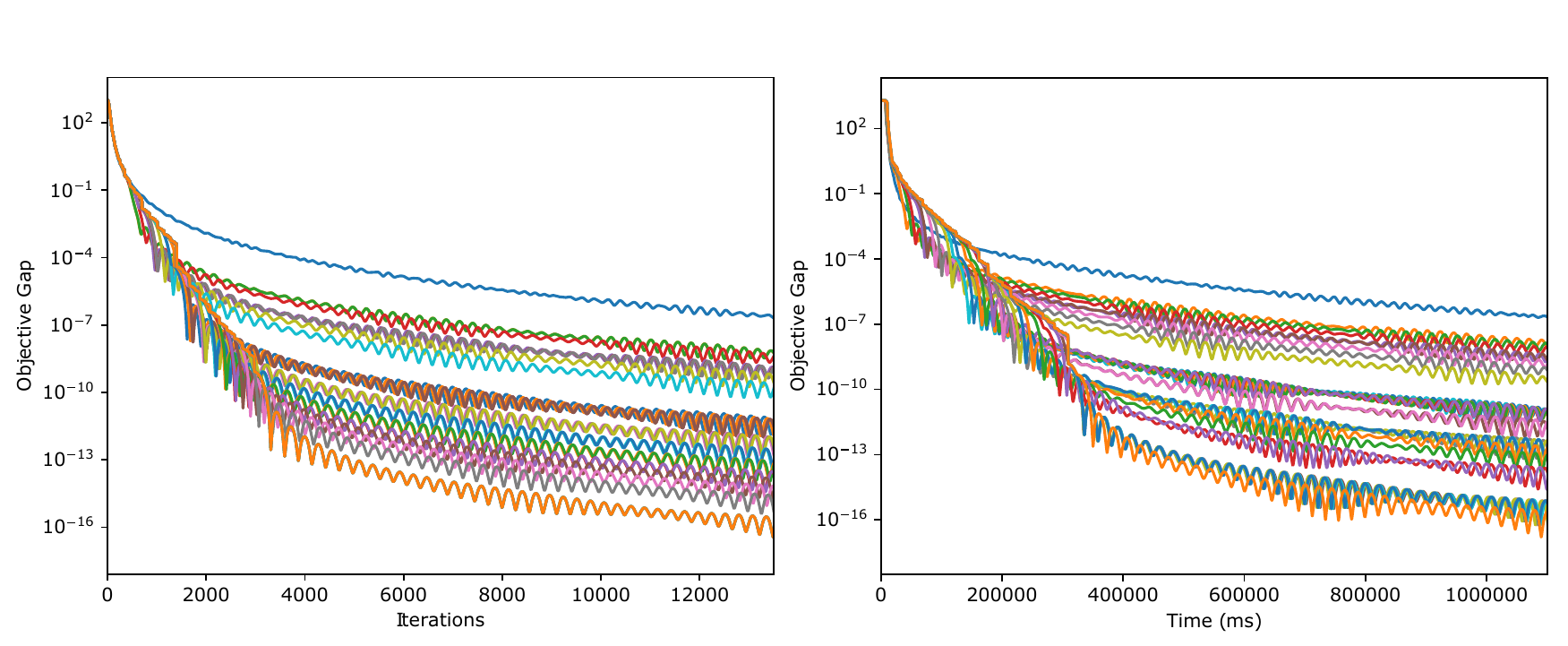}
	\caption{The left plot shows the objective gap of~\eqref{eq:regress} seen by each $\mathtt{accel}$ instance used by $\sparfom$ over 13000 iterations. The right plot shows $\aparfom$ in real time as each $\mathtt{accel}$ instances completes approximately 13000 iterations on an 8-core machine. The highest curve in both plots corresponds to $\mathtt{accel}_{30}$ and the lowest curve corresponds to $\mathtt{accel}_{-1}$ (and everything in between roughly follows this ordering).} \label{fig:accel}
\end{figure}

The results of applying the both the synchronous and asynchronous schemes are shown in Figure~\ref{fig:accel}. Since the $32$ plotted curves heavily overlap and use similar colors, we omit the legend from our plot.
Our implementation of $\aparfom$ launches $32$ threads using the multi-core, distributed processing libraries in Julia (v0.7.0), which each execute one of the schemes $32$ parallel versions of the accelerated method. This code is run on an eight-core Intel i7-6700 CPU, and so multiple $\accel_n$ are indeed able to execute and send messages concurrently. Message passing between processes is done using the RemoteChannel object provided by Julia.

This experiment shows both restarting schemes converge linearly to the target accuracy of $\epsilon=10^{-9}$ within the first $2000$ iterations. Since the first algorithm $\mathtt{accel}_{30}$ never restarts, we see that the accelerated method without restarting only reaches an accuracy of approximately $10^{-3}$ by iteration $2000$.  
Observe that each $\mathtt{accel}_n$ slows down to a sublinear rate after it reaches its target objective gap of $2^n\epsilon$, which corresponds to when that algorithm has restarted for the last time (and will hence not longer benefit from any messages it receives). This behavior matches what our convergence theory predicts (namely, linear convergence to an accuracy of $2^n\epsilon$ followed by the accelerated method's standard $O(\sqrt{L/\epsilon})$ rate).

The results from the asynchronous scheme are fairly similar to those of the synchronous scheme and agree with our theory's predictions. 
We remark that $\aparfom$ is fundamentally nondeterministic since there is a race condition in when each first-order method will complete its assigned task and when messages will be received. Hence, the results in the right side of Figure~\ref{fig:accel} only show one possible outcome of running $\aparfom$, although we found that the general shape of the plot is consistent across many applications of the method.

\bibliographystyle{plain}
\bibliography{a_simple_restart_scheme_v2}

\appendix

\section{{\bf  Complexity of Subgradient Method}}  \label{sect.j}

Here we record the simple proof of the complexity bound for $ \subgrad $ introduced in \S\ref{sect.ba}, namely, that if $ f $ is $ M $-Lipschitz on $ Q $, then for the function 
\begin{equation}  \label{eqn.ja}
  K_{\subgrad}(\delta, \epsilon ) :=  (M \delta/ \epsilon)^2  \; , 
\end{equation}
the iterates $ x_{k+1} = P_Q ( x_k - \frac{\epsilon}{ \| g_k \|^2} g_k ) $ satisfy
\begin{equation}  \label{eqn.jb}
   \dist(x_0, X^*) \leq \delta \quad \Rightarrow \quad  \min \{ f(x_k) \mid 0 \leq k \leq K_{\subgrad}(\delta , \epsilon) \} \, \leq \, f^* + \epsilon \; .   
   \end{equation}

Indeed, letting $ x^* $ be the closest point in $ X^* $ to $ x_k $, and defining $ \alpha_k = \epsilon/\| g_k \|^2 $, we have
\begin{align*} 
     \dist(x_{k+1}, X^*)^2 & \leq \| x_{k+1} - x^* \|^2 
                           = \| P( x_k - \alpha_k g_k) - x^* \|^2 
                      \leq \| (x_k - \alpha_k g_k) - x^* \|^2 \\
                      & = \| x_k - x^* \|^2 - 2 \alpha_k \lin g_k, x_k - x^* \rin + \alpha_k^2 \| g_k \|^2  \\
             & \leq \| x_k - x^* \|^2 - 2 \alpha_k \large(f(x_k) - f^* \large) +  \alpha_k^2 \| g_k \|^2  \quad \textrm{\tiny (because $ g_k $ is a subgradient at $ x_k $)}\\
             & \leq \dist(x_k, X^*)^2 -  \large( \, 2 \large( f(x_k) - f^*   \large)   - \epsilon \, \large) \, \epsilon /M^2 \; ,  
\end{align*}
the final inequality due to $ \alpha_k = \epsilon/ \| g_k \|^2 $, $ \| g_k \| \leq M $ and the definition of $ x^* $. Thus, 
\[  f(x_k) - f^* > \epsilon \,  \Rightarrow \,   \dist(x_{k+1}, X^*)^2 < \dist(x_k, X^*)^2 - (\epsilon/M)^2 \; . \]
Consequently, by induction, if $ x_{k+1} $ is the first iterate satisfying $ f(x_{k+1}) - f^* \leq  \epsilon $, then 
\[ 0 \leq  \dist(x_{k+1},X^*)^2 < \dist(x_0,X^*)^2 - (k+1) (\epsilon/M)^2 \; , \]
  implying the function $ K_{\subgrad} $ defined by (\ref{eqn.ja})  satisfies (\ref{eqn.jb}).
  
  \section{{\bf Smoothing}} \label{sect.k} 
  
  Here we present the first-order method $ \smooth $ which is applicable when the objective function has a smoothing. We also justify the claim that $ K_{\smooth}(\delta, \epsilon) $, as defined by (\ref{eqn.bi}), provides an appropriate bound for the number of iterations. 

We begin with an observation regarding the function $ K_{\fom}(\delta, \epsilon) $ for bounding the number of iterations of $ \fom $ in terms of the desired accuracy $ \epsilon $ and the distance from the initial iterate to $ X^* $, the set of optimal solutions:  Occasionally  in  analyses of first-order methods, the value $ K_{\fom}( \delta_{\epsilon}, \epsilon) $ is shown to be an upper bound, where $ \delta_{\epsilon} $ is the distance from the initial iterate to $ X(\epsilon) := \{ x \in Q \mid f(x) \leq \epsilon \} $. Since $ \delta_{\epsilon} < \delta $, the latter upper bound is tighter.  Particularly notable in focusing on $ \dist(x_0,X(\epsilon)) $ rather than $ X^* $   is the masterpiece \cite{tseng2008accelerated} of Paul Tseng.\footnote{This work can easily be found online, although it was never published (presumably a consequence of Tseng's disappearance).}

In many other cases, even though the analysis is done in terms of $ \dist(x_0,X^*) $, easy modifications can easily be made to obtain an iteration bound in terms of $ \epsilon $ (the desired accuracy) and $ \dist(x_0, X( \epsilon')) $ for $ \epsilon' $ satisfying $ 0 < \epsilon' < \epsilon $. This happens to be the case for the analysis of $ \fista $ in \cite{beck2009fast}  (a special case being $ \accel $, which we rely upon for $ \smooth $. 

In particular, small modifications to the proof of \cite[Thm 4.4]{beck2009fast} show that $ \fista $ computes an $ \epsilon $-optimal solution within a number of iterations not exceeding
\begin{equation}  \label{eqn.ka} 
                      2 \, \dist(x_0,X(\epsilon/2)) \sqrt{L/ \epsilon } \; . 
\end{equation}
We provide specifics in a footnote.\footnote{In the proof of \cite[Thm 4.4]{beck2009fast}, $ x^* $  designates any optimal solution.  However, no properties of optimality are relied upon other than $ f(x^*) = f^* $ (including not being relied upon in \cite[Lemma 4.1]{beck2009fast}, to which the proof refers). Instead choose $ x^* $ to be the point in $ X(\epsilon/2) $ which is nearest to $ x_0 $, and everywhere replace $ f^* $ by $  f^* + \epsilon/ 2 $. 

Also make the trivial observation that the conclusion to \cite[Lemma 4.2]{beck2009fast} remains valid even if the scalars $ \alpha_i $ are not sign restricted, so long as all of the scalars $ \beta_i $ are non-negative (whereas the lemma as stated in the paper starts with the assumption that all of these scalars are positive).     

With these changes, the proof of \cite[Thm 4.4]{beck2009fast}  shows that the iterates $ x_k $ of $ \fista $ satisfy
\[  f(x_k) - (f^* + \epsilon/2)  \leq \frac{2 L \dist(x_0, X(\epsilon/2))^2 }{(k+1)^2} \; . \]
Thus, to obtain an $ \epsilon $-optimal solution, (\ref{eqn.ka}) iterations suffice.}

Recall the following definition:
\begin{quote}
An  ``$ (\alpha, \beta) $ smoothing of $ f $'' is a family of functions $ f_{\eta} $ parameterized by $ \eta > 0 $, where for each $ \eta $, the function $ f_{\eta} $ is smooth on a neighborhood of $ Q $ and satisfies
\begin{enumerate}
\item $ \| \grad f_{\eta}(x) - \grad f_{\eta}(y) \| \leq \smfrac{\alpha }{\eta }  \| x - y \| \; , \quad \forall x,y \in Q $.
\item $ f(x) \leq f_{\eta}(x) \leq f(x) + \beta \eta \; , \quad \forall x \in Q $,
\end{enumerate}
\end{quote}

To make use of an $ (\alpha, \beta) $-smoothing of a nonsmooth objective $ f $ when the goal is to obtain an $ \epsilon $-optimal solution (for $ f $), we rely on $ \accel $ (or more generally, $ \fista $) applied to $ f_{\eta} $ with $ \eta $ chosen appropriately.  In particular, let the algorithm $ \smooth(\epsilon) $ for obtaining an $ \epsilon $-optimal solution of $ f $ be precisely the algorithm $ \accel(2\epsilon/3) $, applied to obtaining an $ (2\epsilon /3) $-optimal solution of $ f_{\epsilon/3\beta} $ (i.e., $ \eta = \epsilon/3\beta $). Note that the Lipschitz constant $ L = \alpha/\eta $ is available for $ \accel $  to use in computing iterates (assuming $ \alpha $ is known, as it is for our examples in \S\ref{sect.ba}).  

To understand these choices of parameters, first observe that by definition, $ \accel(2\epsilon/3 ) $  is guaranteed -- for whatever smooth function to which it is applied -- to compute an iterate which is an $ (2\epsilon/3  ) $-optimal solution.  Let $ x_k $ be such an iterate when the function is $ f_{\epsilon/3 \beta} $, that is  
\begin{equation}  \label{eqn.kb} 
    f_{\epsilon/3 \beta}(x_k) - f^*_{\epsilon /3 \beta} \leq 2 \epsilon /3 \; . 
   \end{equation} 
Now observe that letting $ x^* $ be a minimizer of $ f $, we have
\[  f_{\epsilon/3 \beta}^* \leq f_{\epsilon/3 \beta }(x^*) \leq f(x^*) + \beta ( \epsilon/3 \beta) = f^* + \epsilon/3 \; .   \]
Substituting this and $ f(x_k) \leq f_{\epsilon/3 \beta}(x_k) $ into (\ref{eqn.kb})  gives, after rearranging, $ f(x_k) - f^* \leq \epsilon $, that is, $ \smooth(\epsilon) $ will have computed an $ \epsilon $-optimal solution for $ f $. 

Consequently, the number of iterations required by $ \smooth $ to obtain an $ \epsilon $-optimal solution for $ f $ is no greater than the number of iterations for $ \accel $ to compute an $ (2\epsilon /3) $-optimal solution of $ f_{\epsilon/3 \beta} $, and hence per the discussion giving the iteration bound (\ref{eqn.ka}), does not exceed 
\begin{equation}  \label{eqn.kc} 
2 \, \dist(x_0, X_{\epsilon/3 \beta}( \epsilon/3 ) \sqrt{ (3 \alpha \beta/\epsilon)(2 \epsilon/ 3)}  \; , 
  \end{equation} 
where $ X_{\epsilon/3 \beta}(\epsilon/3) = \{ x \in Q \mid f_{\epsilon/3 \beta}(x) - f^*_{ \epsilon/3 \beta} \leq \epsilon/3 \} $, and where we have used the fact that the Lipschitz constant of $ \grad f_{\epsilon/3 \beta} $ is $ \alpha/ ( \epsilon /3 \beta) $. We claim, however,  $ X^* \subseteq X_{\epsilon/3 \beta}( \epsilon/3)  $, which with (\ref{eqn.kc})  results in the iteration bound
\[   K_{\smooth}(\delta, \epsilon) = 3   \delta \sqrt{2\alpha \beta} / \epsilon 
\] 
for $ \smooth $ to compute an $ \epsilon $-optimal solution for $ f $, where $ \delta $ is the distance from the initial iterate to $ X^* $, the set of optimal solutions for $ f $.

Finally, to verify the claim, letting $ x^* $ be any optimal point for $ f $, we have
\[  f_{\epsilon/3 \beta }(x^*) \leq f(x^*) + \beta (\epsilon/3 \beta) = f^* + \epsilon/3 \leq f_{\epsilon/2 \beta}^* + \epsilon/3 \; , \]
as desired. 

\section{{\bf  Proof of Theorem \ref{thm.ha} } } \label{sect.l}

The proof proceeds through a sequence of results analogous to the sequence in the proof of Theorem~\ref{thm.ea}, although bookkeeping is now  more pronounced. As for the proof of Theorem~\ref{thm.ea}, we refer to $ \fom_n $ ``updating'' at a point $ x $. The starting point $ \mathbf{x_0} $ is considered to be the first update point for every $ \fom_n $.  After $ \aparfom $ has started, then for $ n < N $, updating at $ x $ means the same as restarting at $ x $,  and for $ \fom_{N} $, updating at $ x $  is the same as having computed an iterate $ x $ satisfying the current task of $ \fom_{N} $ (in which case the point is sent to $ \fom_{N-1} $, even though $ \fom_{N} $ does not restart). 

Keep in mind that when a new epoch for $ \fom_n $ occurs, if a message from $ \fom_{n+1} $ arrives in the inbox exactly when the epoch begins, then the restart point will not be decided immediately, due to $ \fom_n $ being made to pause, possibly contiguously. In any case, the restart point is decided after only a finite amount of time, due to $ \fom_{n+1} $ sending at most finitely many messages in total.\footnote{The number of messages sent by $ \fom_{n+1} $ is finite because (1) $ \fom_{n+1} $ sends a message only when it has obtained $ x $ satisfying $ f(x) \leq f(x_{n+1,0}) - 2^{n+1} \epsilon $, where $ x_{n+1,0} $ is the most recent (re)start point, and (2) we assume $ f^* $ is finite.} 

\begin{prop}  \label{prop.la} 
Assume that at time $ t $, $ \fom_n $ updates at $ x $ satisfying $ f(x) - f^* \geq 2 \cdot 2^n \epsilon $. Then $ \fom_n $ updates again, and does so no later than time 
\begin{equation}  \label{eqn.la}
   t + m \, \tpause +  T_{\fom}( D(f(x)), 2^n \epsilon ) \; , 
 \end{equation}   
where $ m $ is the number of messages received by $ \fom_n $ between the two updates.
\end{prop}
\noindent {\bf Proof:}  Let $ \widetilde{m} $ be the total number of messages received by $ \fom_n $ from time zero onward. We know $ \widetilde{m} $ is finite. Consequently, at time
\begin{equation}  \label{eqn.lb} 
   t + \widetilde{m}  \, \tpause +  T_{\fom}( D(f(x)), 2^n \epsilon ) \; ,
\end{equation}    
$ \fom_n $ will have devoted at least $  T_{\fom}( D(f(x)), 2^n \epsilon ) $ units of time to computing iterations.  Thus, if $ \fom_n $ has not updated before time (\ref{eqn.lb}), then at that instant, $ \fom_n $  will not be in a pause, and will have obtained $ \bar{x} $ satisfying  
\begin{equation}  \label{eqn.lc} 
 f( \bar{x}) - f^* \leq 2^n \epsilon \quad \textrm{-- thus, $ f(\bar{x}) \leq f(x) - 2^n \epsilon $ (because $ f(x) - f^* \geq 2 \cdot 2^n \epsilon $)} \; .   
 \end{equation}
Consequently, if $ \fom_n $ has not updated by time (\ref{eqn.lb}), it will update at that time, at $ \bar{x} $.

Having proven that after updating at $ x $, $ \fom_n $ will update again, it remains to prove the next update will occur no later than time (\ref{eqn.la}).   Of course the next update occurs at (or soon after) the start of a new epoch. Let $ m_1 $ be the number of messages received by $ \fom_n $ after updating at $ x $ and before the new epoch begins, and let $ m_2 $ be the number of messages received in the new epoch and before $ \fom_n $ updates (i.e., before the restart point is decided). (Thus, from the beginning of the new epoch until the update, $ \fom_n $ is contiguously paused for an amount of time not exceeding $ m_2 \, \tpause $.) Clearly, $ m_1 + m_2 = m $.  

In view of the preceding observations, to establish the time bound (\ref{eqn.la})  it suffices to show the new epoch begins no later than time
\begin{equation}  \label{eqn.ld} 
   t + m_1 \, \tpause + T_{\fom}( D(f(x)), 2^n \epsilon ) \; . 
 \end{equation}
However, if the new epoch has not begun prior to time (\ref{eqn.ld})  then at that time, $ \fom_n $ is not in a pause, and has spent enough time computing iterations so as to have a point $ \bar{x} $ satisfying (\ref{eqn.lc}), causing a new epoch to begin instantly. \hfill $ \Box $
 
 \begin{prop}  \label{prop.lb}  
 Assume that at time $ t $, $ \fom_n $ updates at $ x $ satisfying $ f(x) - f^* \geq 2 \cdot 2^n \epsilon $. Let $ \mathbf{j}  := \lfloor \frac{ f(x) - f^* }{ 2^n \epsilon } \rfloor - 2 $. Then $ \fom_n $ updates at a point $ \bar{x} $ satisfying $ f( \bar{x}) - f^* < 2 \cdot 2^n \epsilon $ no later than time
\[  
            t + m \, \tpause + \sum_{j=0}^{\mathbf{j}}   T_{\fom}( D( f(x) - j \cdot 2^n \epsilon), 2^n \epsilon) \; ,
\]  
where $ m $ is the total number of messages received by $ \fom_n $ between the updates at $ x $ and  $ \bar{x} $. 
\end{prop} 
\noindent {\bf Proof:} The proof is essentially identical to the proof of Proposition~\ref{prop.fb}, but relying on Proposition~\ref{prop.la}  rather than on Lemma~\ref{lem.fa}.  \hfill $ \Box $

 \begin{cor}  \label{cor.lc} 
Assume that at time $ t $, $ \fom_n $ updates at $ x $ satisfying $ 
 f(x) - f^* < \mathbf{i}  \cdot 2^n \epsilon $,
where $ \mathbf{i}  $ is an integer and $ \mathbf{i} \geq 3 $. Then $ \fom_n $ updates at a point $ \bar{x} $ satisfying  $ f( \bar{x}) - f^* < 2 \cdot 2^n \epsilon $ no later than time
\begin{align}
  & t + m \, \tpause + \sum_{i=3}^{ \mathbf{i}} T_{\fom}(D_{n,i}, 2^n \epsilon) \nonumber \\
   & \qquad  \qquad  \textrm{where} \, \, D_{n,i} = \min \{ D(f^* + i \cdot 2^n \epsilon), D(f(x) ) \} \; ,  \label{eqn.le} 
   \end{align}
and where $ m $ is the number of messages received by $ \fom_n $ between the updates at $ x $ and $ \bar{x} $. 

Moreover, if $ n > -1 $, then after sending the message to $ \fom_{n-1} $ containing the point $ \bar{x} $, $ \fom_n $ will send at most one further message.   
 \end{cor}
 \noindent {\bf Proof:} Except for the final assertion, the proof is essentially identical to the proof of Corollary~\ref{cor.fc}, but relying on Proposition~\ref{prop.lb}  rather than on Proposition~\ref{prop.fb}.
 
For the final assertion, note that if $ \fom_n $ sent two points after sending $ \bar{x} $ - say, first $ x' $ and then $ x'' $ -- we would have
\[   f(x'') \leq f(x') - 2^n \epsilon \leq (f( \bar{x}) - 2^n \epsilon) - 2^n \epsilon < f^* \; , \]
a contradiction.  \hfill $ \Box $

\begin{cor}  \label{cor.ld} 
Let $  n > -1 $. Assume that at time $ t $, $ \fom_{n} $ updates at $ x $ satisfying $ f(x) -  f^* < 5 \cdot 2^n \epsilon $, and assume from time $ t $ onward, $ \fom_n $ receives at most $ \hat{m}   $ messages. Then for either $ \hat{m}'  = 0 $ or $ \hat{m}'  = 1 $, $ \fom_{n-1} $ updates at $ x' $ satisfying $ f( x') - f^* < 5 \cdot 2^{n-1} \epsilon $ no later than time
\[ 
     t + \ttransit + ( \hat{m}  + 2 - \hat{m}' ) \tpause  + \sum_{i=3}^5 T_{\fom}( D_{n,i}, 2^n \epsilon )  
\] 
(where $ D_{n,i} $ is given by (\ref{eqn.le})), and from that time onward, $ \fom_{n-1} $ receives at most $ \hat{m}'  $ messages.
\end{cor}
 \noindent {\bf Proof:}  By Corollary \ref{cor.lc}, no later than time 
\[ 
t + \hat{m} \tpause  + \sum_{i=3}^5 T_{\fom}( D_{n,i}, 2^n \epsilon ) \; , 
\] 
 $ \fom_n $ updates at $ \bar{x} $ satisfying $ f(\bar{x}) - f^* <  2 \cdot 2^n \epsilon $. When $ \fom_n $ updates at $ \bar{x} $, the point is sent to the inbox of $ \fom_{n-1} $, where it arrives no more than $ \ttransit $ time units later, causing a pause of $ \fom_{n-1} $ to begin immediately. Moreover, following the arrival of $ \bar{x} $, at most one additional message will be received by $ \fom_{n-1} $ (per the last assertion of Corollary~\ref{cor.lc}). 
 
If during the pause, $ \fom_{n-1} $ does not receive an additional message, then at the end of the pause, $ \fom_{n-1} $ either updates at $ \bar{x} $ or continues its current epoch.  The decision on whether to update at $ \bar{x} $ is then enacted no later than time
\begin{equation}  \label{eqn.lf} 
 t + \ttransit + (\hat{m}  + 1) \tpause + \sum_{i=3}^5 T_{\fom}( D_{n,i}, 2^n \epsilon ) \; ,
\end{equation} 
after which $ \fom_{n-1} $ receives at most one message.

On the other hand, if during the pause, an additional message is received, then a second pause immediately begins, and $ \bar{x} $ is overwritten by a point $ \bar{\bar{x}} $ satisfying $ f(\bar{\bar{x}}) \leq f( \bar{x}) - 2^n \epsilon  < f^* + 2^n \epsilon $. Here, the decision on whether to update at $ \bar{\bar{x}} $ is enacted no later than time
\begin{equation}  \label{eqn.lg} 
 t + \ttransit + ( \hat{m}  + 2) \tpause + \sum_{i=3}^5 T_{\fom}( D_{n,i}, 2^n \epsilon ) \; ,
\end{equation} 
after which $ \fom_{n-1} $ receives no messages.

If $ \fom_{n-1} $ chooses to update at $ \bar{x} $ (or $ \bar{\bar{x}} $), the corollary follows immediately from (\ref{eqn.lf}) and (\ref{eqn.lg}), as the value of $ f $ at the update point is then less than $ f^* + 4 \cdot 2^{n-1} \epsilon $. On the other hand, if $ \fom_{n-1} $ chooses not to update, it is only because its most recent update point satisfies $ f(x_{n-1,0}) < f( \bar{x}) + 2^{n-1} \epsilon $ (resp., $  f(x_{n-1,0}) < f( \bar{\bar{x}}  ) + 2^{n-1} \epsilon $), and hence $ f(x_{n-1,0}) < f^* + 5 \cdot 2^{n-1} \epsilon $. Thus, here as well, the corollary is established. \hfill $ \Box $
   
 \begin{cor}  \label{cor.le} 
For any $ \bar{N} \in \{ -1, \ldots, N \} $, assume that at time $ t $, $ \fom_{\bar{N}} $ updates at $ x $ satisfying $ f(x) - f^* < 5 \cdot 2^{\bar{N}} \epsilon $, and assume from time $ t $ onward, $ \fom_{\bar{N}} $ receives at most $ \hat{m}  $ messages.  Then no later than time
\[  
     t + (\bar{N} +1) \ttransit + ( \hat{m} + 2 \bar{N} + 2) \tpause    + \sum_{n = -1}^{\bar{N}} \sum_{i=3}^5 T_{\fom}(D_{n,i}, 2^n \epsilon ) \; , 
\] 
$ \aparfom $ has computed an $ \epsilon $-optimal solution.
\end{cor}
\noindent {\bf Proof:}  If $ \bar{N} = -1 $, Corollary~\ref{cor.lc}  implies that after updating at $ x $, the additional time required by $ \fom_{-1} $ to compute a $ (2 \cdot 2^{-1} \epsilon)  $-optimal solution $ \bar{x} $  is bounded from above by
\begin{equation}  \label{eqn.lh} 
    \hat{m} \, \tpause  + \sum_{i=3}^5 T_{\fom}( D_{-1,i}, 2^{-1} \epsilon ) \; . 
    \end{equation} 
The present corollary is thus immediate for the case $ \bar{N} = -1 $.

On the other hand, if $ \bar{N} > - 1 $, induction using Corollary~\ref{cor.ld}  shows that for either $ \hat{m}' = 0 $ or $ \hat{m}' = 1 $,   no later than time
\[     
   t +  (\bar{N}+1) \ttransit + (\hat{m} + 2\bar{N} + 2 - \hat{m}' ) \tpause   + \sum_{n=0}^{\bar{N}} \sum_{i=3}^5 T_{\fom}( D_{n,i}, 2^n \epsilon ) \; , 
\] 
$ \fom_{-1} $ has restarted at $ x $ satisfying $ f(x) - f^* < 5 \cdot 2^{-1} \epsilon $, and from then onwards, $ \fom_{-1} $ receives at most  $ \hat{m}' $ messages.  Then by Corollary~\ref{cor.lc}, the additional time required by $ \fom_{-1} $ to compute a $ (2 \cdot 2^{-1} \epsilon) $-optimal solution does not exceed (\ref{eqn.lh})  with $ \hat{m} $ replaced by $ \hat{m}' $. The present corollary follows.  \hfill $ \Box $
\vspace{2mm}

We are now in position to prove Theorem \ref{thm.ha}, which we restate as a corollary.

\begin{cor}   \label{cor.lf} 
If $ f(\mathbf{x_0}) - f^* < 5 \cdot 2^N \epsilon $, then $ \aparfom $ computes an $ \epsilon $-optimal solution within time
\begin{align} 
  & (\bar{N}+1) \ttransit + 2(\bar{N}+2) \tpause +  3 \sum_{n=-1}^{\bar{N}}  T_{\fom}\left(   D_n , 2^n \epsilon \right)  \label{eqn.li}    \\   
& \qquad \qquad   \textrm{with } \,    D_n :=   \min \{ D(f^* + 5 \cdot 2^n \epsilon), D( f( \mathbf{x_0})) \}  \; , \nonumber 
\end{align}   
where $ \bar{N} $ is the smallest integer satisfying both $ f(\mathbf{x_0}) - f^* < 5 \cdot 2^{\bar{N}} \epsilon $ and $ \bar{N} \geq -1 $.

In any case, $ \aparfom $ computes an $ \epsilon $-optimal solution within time
\[ 
\mathbf{T}_N + 
  T_{\fom}\big( \dist( \mathbf{x_0}, X^*), 2^{N} \epsilon \big)   \; , 
\]  
 where $ \mathbf{T}_N  $ is the quantity obtained by substituting $ N $ for $ \bar{N} $ in (\ref{eqn.li}).
\end{cor}
\noindent {\bf Proof:}  Consider the case $ f( \mathbf{x_0}) - f^* < 5 \cdot 2^N  \epsilon $, and let $ \bar{N} $ be as in the statement of the theorem, that is, the smallest integer satisfying both $ f( \mathbf{x_0}) - f^* < 5 \cdot 2^{\bar{N}}  \epsilon $ and $ \bar{N} \geq -1 $. 

If $ \bar{N} = N $, then $ \fom_{\bar{N}} $ never receives messages, in which case the time bound (\ref{eqn.li})   is immediate from Corollary~\ref{cor.le}  with $ t = 0 $ and $ \hat{m} = 0 $. On the other hand, if $ \bar{N} < N $, then since $ f( \mathbf{x_0}) - f^* < \frac{5}{2}  \cdot 2^{\bar{N}+1} \epsilon $, $ \fom_{\bar{N}+1} $ sends at most two messages from time zero onward. Again (\ref{eqn.li})  is immediate from Corollary~\ref{cor.le}.

Regardless of whether $  f( \mathbf{x_0}) - f^* < 5 \cdot 2^N  \epsilon $, we know $ \fom_{N} $ -- which never restarts -- computes $ x $ satisfying $ f(x) - f^* \leq 2^N \epsilon $ within time 
\begin{equation}  \label{eqn.lj} 
        T_{\fom}( \dist( \mathbf{x_0}, X^*), 2^{N}\epsilon) \; . 
        \end{equation}
If $ x $ is not an update point for $ \fom_{N} $, then the most recent update point $ \hat{x} $ satisfies $ f(\hat{x}) < f(x) + 2^{N} \epsilon \leq  f^* + 2 \cdot 2^N \epsilon $. Hence, irrespective of whether $ x $ is an update point, by the time $ \fom_{N} $ has computed $ x $, it has obtained an update point $ x' $ ($ x' = \hat{x} $ if not $ x' = x $)  satisfying $ f(x') - f^* < 5 \cdot 2^{N} \epsilon $. Consequently, the time required for $ \aparfom $ to compute an $ \epsilon$-optimal solution does not exceed (\ref{eqn.lj})   plus the value (\ref{eqn.li}), where in (\ref{eqn.li}), $ N $ is substituted for $ \bar{N} $.  \hfill $ \Box $

\end{document}